\documentclass[12pt]{article}

\usepackage[utf8]{inputenc}
\usepackage{graphicx}
\usepackage{amssymb}
\usepackage{amsfonts}
\usepackage{amscd}
\usepackage{amsthm}
\usepackage{amsmath}
\usepackage{cite}
\usepackage{float}
\usepackage{url}
\usepackage{times,xcolor}

\graphicspath{ {Figures/} }

\newtheorem{theorem}{Theorem}
\newtheorem{proposition}{Proposition}
\newtheorem{lemma}{Lemma}

\newcommand{\bea}{\begin{eqnarray}}
\newcommand{\eea}{\end{eqnarray}}
\newcommand{\beann}{\begin{eqnarray*}}
\newcommand{\eeann}{\end{eqnarray*}}

\newcommand{\R}{\mathbb{R}}

\newcommand{\eps}{\epsilon}

\title{Dynamics of Non-polar Solutions to the Discrete Painlev\'e I Equation \footnote{This work was partially funded by NSF under grant DMS-1615921 to N. Ercolani.}}

\author
{Nicholas Ercolani, Joceline Lega, Brandon Tippings\\
\normalsize{Department of Mathematics, University of Arizona}}

\begin{document}
\maketitle

\begin{abstract}
This manuscript develops a novel understanding of non-polar solutions of the discrete Painlev\'e I equation (dP1). As the non-autonomous counterpart of an analytically completely integrable difference equation, this system is endowed with a rich dynamical structure. In addition, its non-polar solutions, which grow without bounds as the iteration index $n$ increases, are of particular relevance to other areas of mathematics. We combine theory and asymptotics with high-precision numerical simulations to arrive at the following picture: when extended to include backward iterates, known non-polar solutions of dP1 form a family of heteroclinic connections between two fixed points at infinity. One of these solutions, the Freud orbit of orthogonal polynomial theory, is a singular limit of the other solutions in the family. Near their asymptotic limits, all solutions converge to the Freud orbit, which follows invariant curves of dP1, when written as a 3-D autonomous system, and reaches the point at positive infinity along a center manifold. This description leads to two important results. First, the Freud orbit tracks sequences of period-1 and 2 points of the autonomous counterpart of dP1 for large positive and negative values of $n$, respectively. Second, we identify an elegant method to obtain an asymptotic expansion of the iterates on the Freud orbit for large positive values of $n$. The structure of invariant manifolds emerging from this picture contributes to a deeper understanding of the global analysis of an interesting class of discrete dynamical systems. 
\end{abstract}

\section{Introduction}
 \label{intro}

This article is concerned with the following discrete-time non-auto\-no\-mous dynamical system on the plane $\mathbb{R}^2$:
\begin{align} \label{eq:napdintrotwo} 
x_{n+1}&=\dfrac{n}{Nrx_n}-\dfrac{1}{r}-x_n-y_n \\ 
y_{n+1}&=x_n \nonumber,
\end{align} 
where $r$ and $N$ are parameters. Equation \eqref{eq:napdintrotwo} defines
a sequence of birational maps, $\phi_n$ on the $(x,y)$-plane, for $n \geq 1$, which can be extended to $n \leq 0$.  For $n \ne 0$ these maps are undefined along the $y$-axis; however, this does not affect our principal considerations, as we will explain.
In many areas of mathematics \eqref{eq:napdintrotwo} is well known as the discrete Painlev\'e I equation (or just dP1).  There is a vast literature on the relevance of Painlev\'e equations, both continuous and discrete, to various subfields of mathematics and we refer the reader to \cite{bib:cm20}, \cite{bib:kny17}, and \cite{bib:va18} for recent reviews of many of these connections. 

Painlev\'e equations in general were originally characterized by the special, restricted behavior of the singularities that their solutions may have. Indeed, in the autonomous limit of  (\ref{eq:napdintrotwo}), which we may specify by setting $N = \alpha^{-1} \, n$ for some non-zero real constant $\alpha$, the system is analytically completely integrable, by which we mean that its solutions lie on closed curves and may be explicitly written in terms of meromorphic functions. For instance, in the case of autonomous dP1, solutions are expressed in terms of rational functions on an elliptic curve or one of its degenerations \cite{bib:be19}. For general continuous Painlev\'e equations, this translates to solutions being expressible in terms of functions whose global analytic continuations are restricted only to have poles as singularities. This criterion is known as \textit{the Painlev\'e property}. For general discrete Painlev\'e equations, this property can be re-expressed dynamically in terms of {\it singularity confinement}: solutions may become arbitrarily large in a finite number of steps before returning to a bounded region of the plane \cite{bib:grp91,bib:lg03}. That number of steps is fixed for almost all solutions and excursions to infinity can only take place along a fixed set of directions in the plane. 

System \eqref{eq:napdintrotwo} is not known to be analytically completely integrable and the existence of a closed form analytic expression of its solutions remains an open problem. Nevertheless, dP1 continues to possess many if not all of the well-known properties typically associated with integrable systems, such as singularity confinement of its generic solutions \cite{bib:grp91,bib:otgr99,bib:gr04}, zero algebraic entropy \cite{bib:bv99,bib:otgr99}, as well as Lax pair \cite{bib:fik91,bib:pap92} and Hirota birational \cite{bib:kny17} formulations. We will refer to solutions that leave the first quadrant and exhibit singularity confinement as {\it polar solutions}. There are, however, certain initial conditions leading to solutions that increase without bound for all time.  We will refer to these as {\it non-polar} solutions, similar to the \textit{pole-free} terminology originally used in \cite{bib:jos97}. We do not claim that this dichotomy is exhaustive although numerical explorations (a few examples of which are provided in this manuscript) suggest this is the case. The primary focus of this paper is on non-polar solutions.

An important example of a non-polar solution, which originally motivated our interest in this work, stems from the connection of dP1 to approximation theory, more explicitly to analyzing the structure of orthogonal polynomials with exponential weights. 
Recall that a family of orthonormal polynomials $\{ p_{m}\}$ associated to an exponential weight 
$w(\lambda) = \exp(-V(\lambda))$ is a complete basis of polynomials for the weighted $L^2$ space such that 
\begin{eqnarray*}
\int_{\mathbb{R}} p_n(\lambda) p_m(\lambda) w(\lambda) d \lambda &=& \delta_{nm}.
\end{eqnarray*}
It is known that such polynomials are algebraically specified as solutions to a three-term recurrence relation 
with real coefficients. When $V$ is even, this recurrence takes the form
\bea \label{eq:threetermintro}
\lambda p_n(\lambda) = b_{n+1} p_{n+1}(\lambda) + b_n p_{n-1}(\lambda)
\eea
and the recurrence coefficients $\{b_n\}$ themselves are solutions to a 
{\it nonlinear} difference equation (see for instance \cite{bib:va18}). In the particular case where $V(\lambda) = N \left( \frac12 \lambda^2 + \frac{r}4 \lambda^4\right)$ \cite{bib:fre76,bib:mag99,bib:va18}, this difference equation reads
\begin{equation}
\label{eq:Freudeq}
rb_n^2\left( b_{n+1}^2 + b_n^2 + b_{n-1}^2 \right) + b_n^2 = \frac{n}{N},
\end{equation}
which, setting $x_n = b_n^2$, is precisely equivalent to the dP1 system (\ref{eq:napdintrotwo}). Freud \cite{bib:fre76} went on to pin down that $x_n \propto n^{1/2}$, using properties of orthogonal polynomials. We call the corresponding orbit $\{x_n = b_n^2\}$, the \textit{Freud orbit}.
Subsequently, M\'at\'e, Nevai and Zaslavsky \cite{bib:mnz85} built on Freud's work to show the existence of a full asymptotic expansion of the orthogonal polynomial solution in powers of $n^{1/2}$, although they did not describe the coefficients in that expansion. Later Ercolani, McLaughlin, and Pierce \cite{bib:emp08} developed an alternative approach to the asymptotic expansion of the recursion coefficients using Riemann-Hilbert analysis. This work provided a geometric characterization of the recursion coefficients, in terms of a graphical enumeration problem related to diagrammatic expansions in mathematical physics. Since $x_n = b_n^2 >0$, Freud's special orbit remains in the first quadrant for all time. This feature motivated Lew and Quarles \cite{bib:lq83} to seek {\it all} solutions of
\begin{equation}
\label{eq:LQN}
x_n\, (x_{n+1} + x_n + y_n) = n, \quad y_{n+1} = x_n, \quad n > 0
\end{equation}
that possess this property. For
$x_n \ne 0$, \eqref{eq:LQN} is equivalent to \eqref{eq:napdintrotwo} without the term in $1/r$ on the right-hand side of the equation for $x_{n+1}$. In \cite{bib:lq83}, Lew \& Quarles established, by an elegant contraction mapping argument, that there is a one parameter family of such non-polar solutions. According to \cite{bib:ansva15}, addressing the existence and uniqueness of the positive solution of \eqref{eq:LQN} with $y_1 = 0$ was a problem posed by Nevai, and independently solved by him in \cite{bib:n83}. In Appendix \ref{app:LQ_construction}, we extend the Lew-Quarles construction to describe a broad class of non-polar orbits of dP1, which we refer to as the Lew-Quarles orbits.

\begin{figure}[ht]
\centering
\includegraphics[width=\textwidth]{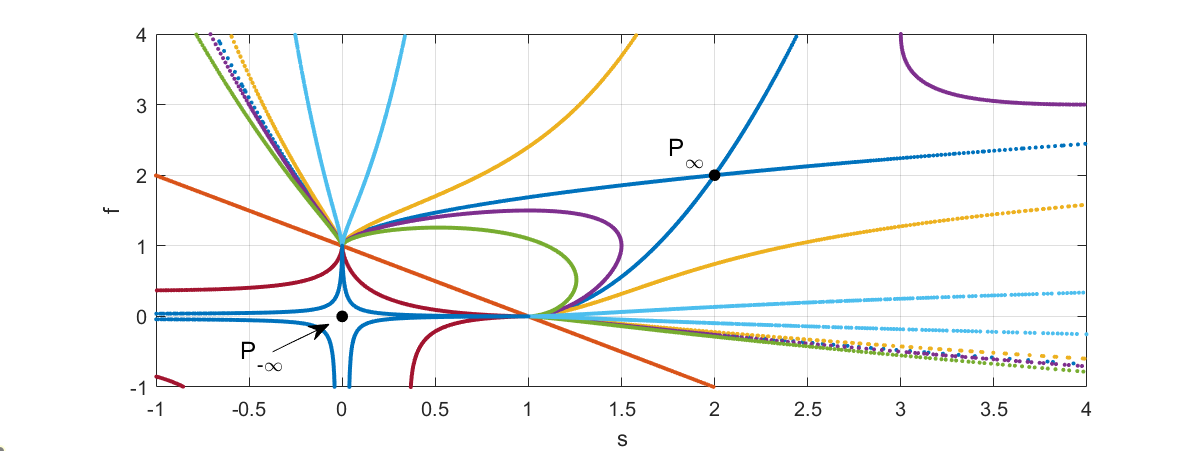}
\caption{Phase portrait of the discrete system \eqref{eq:u=0_plane} defined in the invariant plane at infinity $u = 0$. Different colors correspond to different orbits (most of which repeatedly escape the field of view), each with $3000$ iterates.}
\label{fig:u=0}
\end{figure} 

In this paper, we adopt a dynamical systems perspective for \eqref{eq:napdintrotwo}. Doing so enables us to simultaneously describe and compare polar and non-polar solutions, thereby illuminating novel features of the fuller class of non-polar solutions of dP1. To this end, we use the Painlev\'e property to complete the phase space at infinity with an asymptotic change of variables based on the Riemann-Hilbert scalings used in \cite{bib:emp08}. A further transformation allows us to define the
{\sl asymptotic change of coordinates} $(x,y,n) \to (s, f, u)$, 
\begin{equation}
\label{eq:xyn_usf}
s = \frac{y}{x} + 1 + \frac{1}{r x}, \qquad f = \frac{n}{N r x^2} - \frac{y}{x}, \qquad u = - \frac{1}{r x},
\end{equation}
which reveals the existence of an {\it invariant plane at infinity} ($u=0$). The dynamics in this invariant plane, which organizes the various asymptotic behaviors of the full system, is shown in Figure \ref{fig:u=0} and corresponds to the reduced discrete system 
\begin{equation}
\label{eq:u=0_plane}
s_{n+1} = Z_n f_n, \quad
f_{n+1} = Z_n^2 s_n, \quad Z_n = (f_n -1)^{-1}.
\end{equation}
The picture we will develop is that the Freud orbit, when extended to negative discrete times, is a singular heteroclinic connection from the origin (which we call $P_{- \infty}$) of the $(s,f,u)$ space to the fixed point $P_\infty$ with coordinates $s = f = 2$ and $u=0$. In addition, the Lew-Quarles orbits form a family of non-singular heteroclinic connections that converge to the Freud orbit, both as $n \to -\infty$ and $n \to \infty$. The asymptotic points $P_{-\infty}$ and $P_\infty$ lie in the plane $u=0$ and are marked as black dots in Figure \ref{fig:u=0}. Given the phase plane structure suggested by this figure, heteroclinic connections from $P_{-\infty}$ to $P_\infty$ are trajectories that ``escape'' into the third dimension before asymptotically returning to the invariant plane. Our goal is to characterize the structure of these connections. 

The rest of this article is organized as follows. Section \ref{sec:FO} defines the Freud orbit, explains how it is initialized, and presents a high precision numerical simulation of forward and backward iterates from the selected initial condition, illustrating convergence to $P_{-\infty}$ as $n \to -\infty$, and $P_\infty$ as $n \to \infty$. Section \ref{sec:usf_coord} introduces the change of variables that transforms system \eqref{eq:napdintrotwo} into an autonomous system in $(s,f,u)$ coordinates, and reviews the basic properties of the resulting discrete dynamical system. Section \ref{sec:NP} presents a detailed study, both analytical and numerical, of the Lew-Quarles orbits, and characterizes these solutions as heteroclinic connections between $P_{-\infty}$ and $P_\infty$. In particular, we provide strong evidence that these trajectories live on the center and center-stable manifolds of these two points respectively, and converge exponentially to the Freud orbit as $n \to \infty$. We also build on these results to propose a unique characterization of the Freud orbit as a singular limit of the Lew-Quarles orbits. In Section \ref{sec:Freud}, we prove that orbits that converge to $P_{-\infty}$ and $P_\infty$ along specific invariant curves track sequences of points formed by the period-2 points (as $n \to -\infty)$ and fixed points (as $n \to \infty$) of the autonomous dP1 system, in which $n$ appears as a parameter. This system is further described in Appendix \ref{app:Aut}. In addition, we obtain expansions of these orbits in powers of $|n|^{-1/2}$, valid as $n \to -\infty$ and $n \to \infty$. Applying these results to the Freud orbit, as legitimized by our numerical observations, provides a simple means to obtain its asymptotic expansions to arbitrary order in powers of $|n|^{-1/2}$ as $n \to -\infty$ or $n \to \infty$. Finally, Section \ref{sec:conclusion} summarizes our findings and identifies future directions that build on the dynamical perspective introduced in this work. Numerical methods are presented in Appendix \ref{app:Num}, and asymptotic expansions for the invariant curves near $P_{-\infty}$ and $P_\infty$ are given in Appendix \ref{app:CM}.

\section{Definition of the Freud orbit}
\label{sec:FO}

\begin{figure}[hbtp]
\centering
\includegraphics[width=\textwidth]{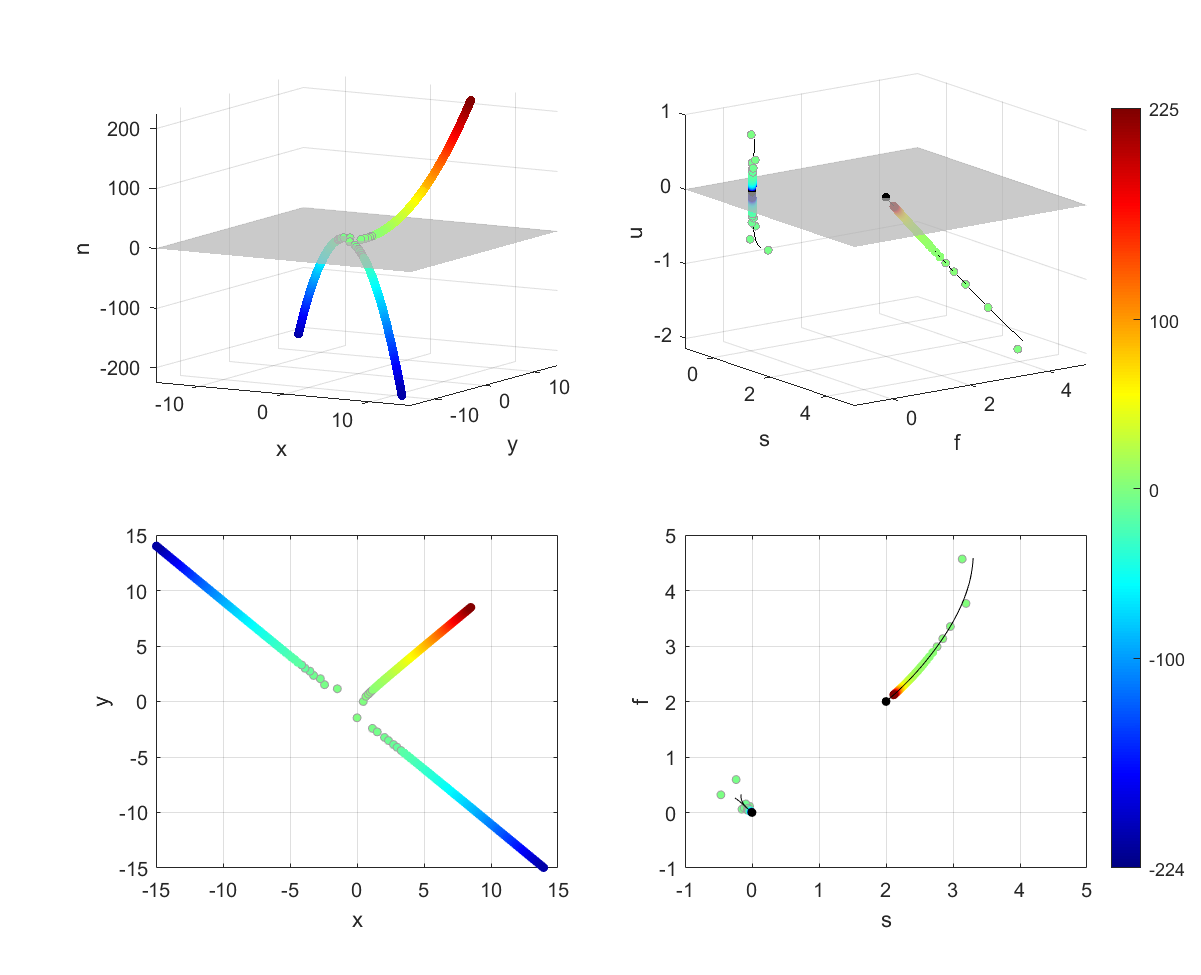}
\caption{Forward and backward iterates of the Freud orbit with $r = N = 1$ in the $(x,y,n)$ coordinates (top left panel) and in the asymptotic coordinate system $(s,f,u)$ (top right panel). Values of $n$ range from $-224$ to $225$. The bottom panels are projections on the $(x,y)$ (left) and $(s,f)$ (right) planes (highlighted in gray in the top panels). Colors range from dark blue to dark red as $n$ increases (see color bar). For negative values of $n$, iterates alternate between the two branches (in blue) in the second and fourth quadrants of the $(x,y)$ plane. Points corresponding to integer values of $n$ near $0$ are in light green and are contoured in grey for added visibility. In the right column, the black dots represent $P_{-\infty}$ and $P_\infty$ (top), as well as their projections (bottom). A 10\textsuperscript{th} order expansion (Equations \eqref{eq:CMP2}) of an invariant curve transverse to the plane $u=0$ valid near $P_{-\infty}$, and a 6\textsuperscript{th} order expansion (Equations \eqref{eq:CMP1}) valid near $P_\infty$, are also plotted (solid black curves). The reader is referred to Appendix \ref{app:Num} for details on the numerical simulations.}
\label{fig:comb}
\end{figure}

As mentioned above, Freud's asymptotic analysis of the recurrence coefficients for orthogonal polynomials with quartic weight (\ref{eq:threetermintro}) determines a particular solution of (\ref{eq:napdintrotwo}) with the property that in forward time the orbit always lies in the first quadrant of the phase plane and, further, that $x_n \propto \sqrt{n}$. Using the orthonormality of the $p_n$ (see Appendix \ref{app:OP}), one finds that the initial conditions for this orbit are  
\begin{equation}
y_1^F = b_0^2 = 0, \ \ \ x_1^F=b_1^2=\dfrac{\mu_2}{\mu_0}, \label{eq:introxnintialcondition}
\end{equation} where the $\mu_i$ are the moments:
\begin{equation}
\mu_i=\int_{\mathbb{R}} \lambda^i w(\lambda) d\lambda. \label{eq:intromomentsdef}
\end{equation}
We will refer to  iterates of $(x_1^F, y_1^F)$ under \eqref{eq:napdintrotwo} as the {\it Freud orbit} and use the superscript $F$ to distinguish this specific set of points. Another form of initial conditions is given by $(x_2^F, x_1^F)$ (See Appendix \ref{app:OP}).

Figure \ref{fig:comb} shows a high-precision numerical simulation of the Freud orbit in the $(x,y,n)$ (top left) and $(s,f,u)$ (top right) coordinates, for $n \in [-224, 225]$ and parameter values $r = N = 1$. The color scale goes from dark blue (large, negative values of $n$) to dark red (large, positive values of $n$), with lighter colors corresponding to negative and positive values of $n$ that are smaller in magnitude. The bottom row shows projections of the orbit on the $(x,y)$ plane (left), and on the $(s,f)$ plane (right). In the right column, the points, $P_{-\infty}=(0,0,0)$ and $P_\infty = (2,2,0)$, as well as their projections on the $(s,f)$ plane, are marked with black dots. In the bottom left panel, the point $(x_1^F,y_1^F)$ is shown in bright green on the $x$-axis with $x_1^F \simeq 0.47$. The image of $(x_1^F,y_1^F,n=1)$ in the $(s,f,u)$ space (top right panel) is the point of approximate coordinates $(3.1,4.6,-2.1)$ in the bottom right corner of the plot. As detailed in Section \ref{sec:NP}, we will define the pre-image of $(x_1^F,y_1^F)$ as the point on the $y$-axis of coordinates $x_0^F = 0$, $y_0^F \simeq -1.5$. In the $(s,f,u)$ space, iterates blow up when $n=0$ but are well defined for $n \ne 0.$ The forward part ($n > 0$) of Freud's orbit is in green, yellow, and red, with $x \simeq y \simeq n^{1/2}$. Its backward iterates (negative values of $n$) alternate between the second and fourth quadrants (points shown in green as well as light and dark blue). In the asymptotic coordinate system (right column), the Freud orbit moves away from $P_{-\infty}$ while alternating between each side of the plane $u=0$, and converges to $P_\infty$ as $n \to \infty$. The thin solid curves correspond to Equations  \eqref{eq:CMP1} and  \eqref{eq:CMP2} (see Section \ref{sec:NP}), which capture the dynamics near $P_\infty$ and $P_{-\infty}$ respectively. 

As is clearly illustrated in Figure \ref{fig:comb}, the $(s,f,u)$ coordinate system provides a natural framework to analyze the properties of the Freud orbit. We now explain the origins of this change of coordinates.

\section{Asymptotic change of coordinates}
\label{sec:usf_coord} 


The non-autonomous dP1 mapping \eqref{eq:napdintrotwo} may be transformed into an autonomous 3-dimensional system by introducing the variable $\alpha_n = n / N$. The choice of denominator here is motivated by a fundamental re-scaling used in the Riemann-Hilbert analysis of \cite{bib:emp08}, which amounts to considering a limit in which $n$ and $N$ go to infinity together at a fixed rate. We first set
\[
\theta_1 = \psi^2 \alpha, \qquad \theta_2 = \psi \sqrt N y,  \qquad \psi = \frac{1}{\sqrt N x}.
\]
The singling out of the coordinates $\theta_1 = n / (N x)^2$ and $\theta_2 = y / x$ stems from an alternative approach to extending the phase space at infinity developed in \cite{bib:kny17}, based on a resolution of singularities (see also \cite{bib:tip20}). The above change of variables leads to
\begin{align}
\theta_{1,n+1} &= Z_n^2 \left(\theta_{1,n} + \frac{1}{N} \psi_n^2\right) \nonumber \\
\theta_{2,n+1} &= Z_n = \gamma \left(\theta_{1,n} - \frac{1}{\sqrt N} \psi_n - \gamma - \gamma \theta_{2,n}\right)^{-1} \label{eq:infty} \\
\psi_{n+1} &= Z_n \psi_n, 
\nonumber
\end{align}
where $\gamma = r/N.$ This system has two fixed points, 
\[
P_\infty = (3 \gamma, 1, 0), \qquad P_{-\infty} = (-\gamma, -1, 0),
\]
corresponding to positive ($P_\infty; \ \theta_1 > 0$) and negative ($P_{-\infty}; \ \theta_1 < 0$) values of $n$. The associated values of $Z$ are $1$ and $-1$, respectively. The eigenvalues and eigenvectors of the linearization of \eqref{eq:infty} about these fixed points are
\begin{align*}
P_\infty:\quad & \lambda = 1, & e_1 &= \left(1, 0, \sqrt N\right)^T \\
& \lambda = -2 \pm \sqrt 3, & e_\pm &= \left(\gamma (3 \mp \sqrt 3), 1, 0\right)^T \\
P_{-\infty}:\quad & \lambda = -1, & \xi_1 &= \left(\gamma, 1, -\gamma \sqrt N\right)^T\\
& \lambda = \pm i, & \xi_\pm &= \left(\gamma (1 \mp i), 1, 0\right)^T.
\end{align*}
We now introduce a change of coordinates centered on the fixed point $P_{-\infty}$ and consistent with the basis of eigenvectors of the linearization of \eqref{eq:infty} near $P_{-\infty}.$ Specifically, setting
\begin{equation}
\theta_1 = - \gamma + \gamma (u + s + f), \quad \theta_2 = - 1 + u + s, \quad \psi = - \gamma \sqrt N u,
\label{eq:usf_coord}
\end{equation}
transforms \eqref{eq:infty} into the following discrete dynamical system
\begin{align}
s_{n+1} &= Z_n f_n, \qquad Z_n = (u_n + f_n -1)^{-1} \nonumber\\
f_{n+1} &= Z_n^2 \left(s_n + \gamma u_n^2 \right), \label{eq:usf_sys}\\
u_{n+1} &= Z_n u_n \nonumber
\end{align}
with which we will now work.

System \eqref{eq:usf_sys} has exactly two fixed points, which in the $(s,f,u)$ coordinates are given by
\[
P_{-\infty} = (0, 0, 0),  \qquad P_\infty = (2, 2,0),
\]
and therefore lie in the invariant plane $u = 0$. Looking for periodic orbits of higher order, we find by direct calculation that there is only one genuine period-2 orbit, given by
\[
(2, 0, 0) \to (0, 2, 0) \to (2, 0, 0),
\]
and a line of period-3 orbits such that $u_n = 0$ and $s_n + f_n = 1$, with $s_n \ne 0$ and $f_n \ne 0$. They are of the form
\begin{align}
\label{eq:P3} \left(s_0, 1 - s_0, 0 \right) \to \left(1 - s_0^{-1}, s_0^{-1}, 0\right) & \to \left((1 - s_0)^{-1}, \left(1 - s_0^{-1}\right)^{-1}, 0\right) \\ & \to \left(s_0, 1 - s_0, 0 \right), \nonumber
\end{align}
where $s_0 \in \R \setminus \{0,1\}.$ Moreover, these are the only period-3 orbits in the invariant plane $u = 0$, other than the two fixed points. The points $s_n = 0$ and $f_n = 0$ are special, in the sense that they are part of a singular period-3 orbit of the form
\begin{equation}
\label{eq:so}
(0, 1, 0) \to (\infty, - \infty, 0) \to (1, 0, 0) \to (0, 1, 0).
\end{equation}
Our numerical explorations, in which we set $\gamma=1$, indicate that polar orbits (see examples in Appendix \ref{app:Num}) track the phase plane structure shown in Figure \ref{fig:u=0} as soon as $|x|$ is larger than a few units (thus corresponding to $|u| < 0.5$). We also observe almost periodic orbits associated with large values of $x$ and $y$ ($|x|,\ |y| > 1000$), which in $(s,f,u)$ coordinates correspond to solutions close to \eqref{eq:P3} on the line $s + f = 1$; in that case, $u$ remains small ($|u| < 10^{-3}$) and a slight drift is observed across the line $s+f = 1$. Polar orbits for which there exists a value of $n$ such that $x_n = {\mathcal  O}(\epsilon)$ and $y_n = {\mathcal  O}(1)$ display generic singularity confinement: iterates of $x_n$ become large before returning to values of order $y_n$. Specifically,
\begin{align*}
x_{n-1} = y_n &= {\mathcal  O}(1), \quad x_n = \epsilon, \quad x_{n+1}=\frac{n}{N r \epsilon}+{\mathcal  O}(1), \quad x_{n+2}=-\frac{n}{N r \epsilon} + {\mathcal  O}(1), \\ 
&x_{n+3} = {\mathcal  O}(\epsilon), \quad x_{n+4} = {\mathcal  O}(y_n) = {\mathcal  O}(1), \quad \text{as } \epsilon \to 0.
\end{align*}
In $(s,f,u)$ coordinates, this translates to
\begin{align*}
&s_{n+1} = 1 + \frac{N \epsilon}{n} + {\mathcal O}(\epsilon^2), \quad f_{n+1} = {\mathcal O}(\epsilon^2), \quad u_{n+1} = - \frac{N \epsilon}{n} + {\mathcal  O}(\epsilon^2)\\
&s_{n+2}= {\mathcal O}(\epsilon^2), \quad f_{n+2} = 1 - \frac{N \epsilon}{n} + {\mathcal O}(\epsilon^2), \quad u_{n+2} = \frac{N \epsilon}{n} + {\mathcal  O}(\epsilon^2),
\end{align*}
as $\epsilon \to 0$. In other words, during singularity confinement, iterates of dP1 visit neighborhoods of the points $(1,0,0)$ and $(0,1,0)$ of the singular period-3 orbit \eqref{eq:so}.

The linearization about the line of period-3 orbits \eqref{eq:P3} is given by 
\begin{align*}
s_{n+3}-s_n &= - 2 \mu_n - 3 \nu_n  -3 u_n + O(\eps) \\
f_{n+3} - f_n &= 3 \mu_n + 4 \nu_n + 4 u_n + O(\eps) \\
u_{n+3}-u_n &= O(\eps)
\end{align*}
for $s_n = s_0 + \mu_n$, $f_n = 1 - s_0 + \nu_n$, $u_n = O(\eps)$, $\mu_n = O(\eps)$, and $\nu_n = O(\eps).$
The matrix
\[
A = \left(\begin{array}{ccc} -2 & -3 & -3 \\ 3 & 4 & 4 \\ 0 & 0 & 0 \end{array}\right)
\]
has eigenvalues $\lambda_A = 0$ with  eigenvector $\zeta_0 = (0, -1,1)^T$ and $\lambda_A = 1$ with eigenvector $\zeta_1 = (1, -1,0)^T$. In addition, since $\lambda_A = 1$ has algebraic multiplicity 2 but geometric multiplicity 1, we define the generalized eigenvector $\zeta_2$  such that 
\[
\zeta_2 = (-1/3,0,0), \quad (A - I) \zeta_2 = \zeta_1.
\]
With these conventions, any initial condition of the form $\zeta = a\, \zeta_0 + b\, \zeta_1 + c\, \zeta_2$ evolves according to
\[
A^k \zeta = \left(\begin{array}{c} b + k c - \frac{c}{3} \\ - (b + k c) \\ 0 \end{array}\right),
\]
so that a perturbation transverse to the line of period-3 orbits, of the form $\zeta = a\, \zeta_0 + c\, \zeta_2$ will linearly drift along this line according to $s_{3 p} = s_0 + 3\, p\, c - c/3$, $f_{3 p} = 1 - s_0 - 3\, p\, c$, $u_{3p} = 0$, and eventually reach the vicinity of the singular orbit \eqref{eq:so}.

Finally, we note that the change of coordinates from $(s,f,u)$ to $(x,y,n)$ is given by
\begin{equation}
\label{eq:usf_xyn}
x=-\frac{1}{r u}, \qquad y = - \frac{s+u-1}{r u}, \qquad \alpha = \frac{s+f+u-1}{r u^2},
\end{equation}
with $\alpha = n/N$.

\section{Characterization of the Lew-Quarles orbits}
\label{sec:NP}
\begin{figure}[ht]
\centering
\includegraphics[width=.8\textwidth]{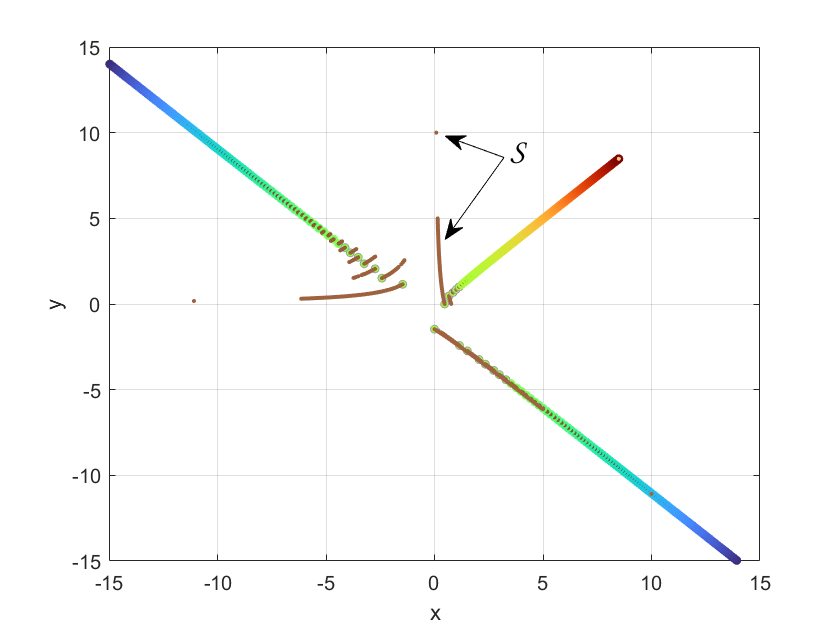}
\caption{Forward and backward iterates of the Lew-Quarles orbits in the $(x,y)$ plane. A subset of the set $\mathcal S$ of initial conditions is shown in brown (along the curve with $y>0$ and $0<x<x_1^F \simeq 0.47$). Iterates of $\mathcal S$ under \eqref{eq:napdintrotwo} (also in brown) remain in the first quadrant and quickly converge to the Freud orbit (in green, orange, and red). The pre-image of $\mathcal S$, or equivalently $\psi_0({\mathcal S})$, is the set of brown colored points in the fourth quadrant. Images of $\mathcal S$ under $\psi_n,\ n \le -1$ are in brown and alternate between the second and fourth quadrants. They are only visible in the second quadrant because they are are all ``aligned'' with, and thus hidden by, $\psi_0({\mathcal S})$ in the fourth quadrant. Backward iterates of $\mathcal S$ converge to the backward iterates of the Freud orbit, shown in green and blue.}
\label{fig:xyLQ}
\end{figure} 
We now restrict our attention to the family of Lew-Quarles orbits, whose proof of existence is summarized in Appendix \ref{app:LQ_construction}. The goal of this section is to show that  these orbits, which are initiated in the first quadrant, can be extended to heteroclinic connections between $P_{-\infty}$ and $P_\infty$, and that the Freud orbit is a singular limit of such solutions. For reference, Figure \ref{fig:xyLQ} shows the Lew-Quarles (in brown) and Freud (in color) orbits in the $(x,y)$ plane.

\subsection{Forward iterates}
In Appendix \ref{app:LQ_construction} we provide an adapted version of Lew and Quarles'  construction that leads to the following existence and uniqueness theorem for a broad class of non-polar orbits for dP1. A different proof, for a more general class of systems that includes dP1, is provided in \cite{bib:ansva15}.
\begin{theorem}
For any $\xi_0 = y_1  \geq 0$ there is a unique solution of (\ref{eq:napdintrotwo}) that remains in the first quadrant. It is defined in terms of a sequence $\{\xi_n\}_{n \ge 0}$ which is a fixed point of a contraction mapping, has initial value $(x_1, y_1) = (\xi_1, \xi_0)$, and is such that 
$\xi_n = x_n = y_{n+1} > 0$ for all $n \ge 0$. 
\label{th:LQ}
\end{theorem}
\noindent We define $\mathcal{S}$ to be the set of initial conditions described in Theorem \ref{th:LQ}. As explained in Appendix \ref{app:Num}, any point on $\mathcal S$, including the initial condition for the Freud orbit $F$, may be numerically estimated with arbitrary precision using the Lew-Quarles contraction mapping. Then, we repeatedly apply the forward mapping $\phi_n, \ n\ge 1$ defined in \eqref{eq:napdintrotwo}, to find the associated orbit. Figure \ref{fig:xyLQ} shows that as $n$ increases, the Lew-Quarles orbits quickly collapse onto $F$. It is natural to conjecture that such a convergence is exponentially fast. This is confirmed by our simulations, which in fact provide strong numerical evidence that these orbits approach $F$ at a {\it uniform} exponential rate. 

\begin{figure}[ht]
\centering
\includegraphics[width=.9\textwidth]{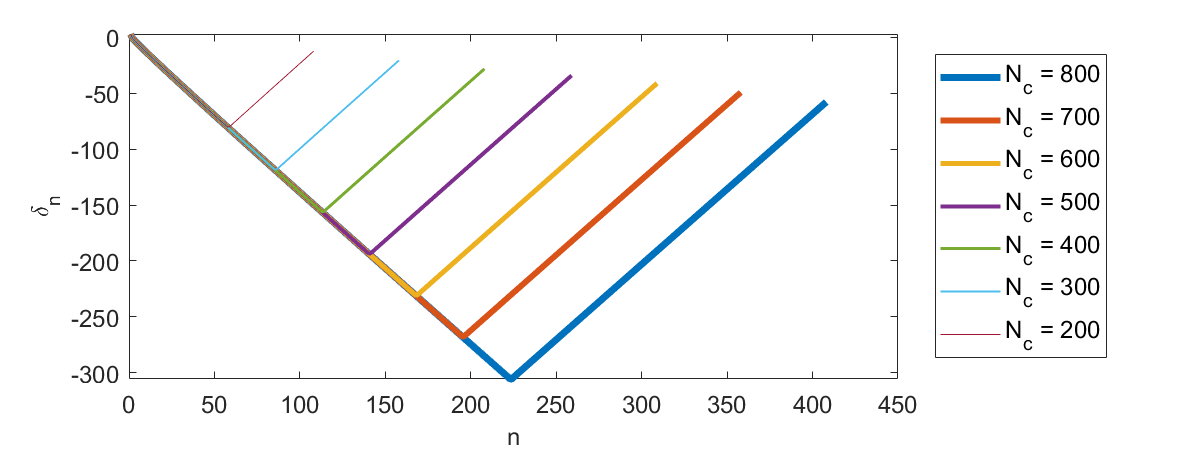}
\caption{Log distance $\delta_n$ between the orbit of a point $Q$ on $\mathcal{S}$ (corresponding to $\xi_0 = 20$) and the Freud orbit, for varying number of contractions $N_c$ used to numerically compute the initial condition 
$Q \in \mathcal{S}$.}
\label{fig:turnaround}
\end{figure} 

Figure \ref{fig:turnaround} shows how the Lew-Quarles orbit originating from a numerically approximated initial condition $Q$ on $\mathcal{S}$, in this case prescribed by $\xi_{0} = 20$, converges exponentially to $F$ in forward time. The $y$-axis displays the quantity:
\begin{equation}
\delta_n = \log || \Phi^n (Q)- \Phi^n(x_1^F,y_1^F) ||,
\end{equation} 
the log distance between the $n\textsuperscript{th}$ iterate of the orbit of $Q$ and the $n\textsuperscript{th}$ iterate of the Freud initial condition ($\xi_0=0$), as a function of $n$ (the iterate number). With our previously introduced notation, $\Phi^n = \phi_n \circ \cdots \circ \phi_1$. Different traces, with increasing numbers of Lew-Quarles contractions used to define the point $Q$, are provided to illustrate how $\delta_n$ behaves under improved initial condition calculation. In all of these instances, the initial condition for the Freud orbit is calculated with 800 Lew-Quarles contractions. Turnaround in the trend, at the vertices of each of the graphs shown in Figure \ref{fig:turnaround}, is indicative of accumulation of numerical error, as this turnaround time increases with the number of contractions used to estimate $Q$. These graphs are truncated just before the approximated orbit exits the first quadrant, indicating a critical accumulation of numerical error (the true orbits never exit the first quadrant, by virtue of $Q$ being on $\mathcal{S}$, the set of initial conditions that, by Theorem \ref{th:LQ}, lead to positive iterates).

\begin{figure}[ht]
\centering
\includegraphics[width=.9\textwidth]{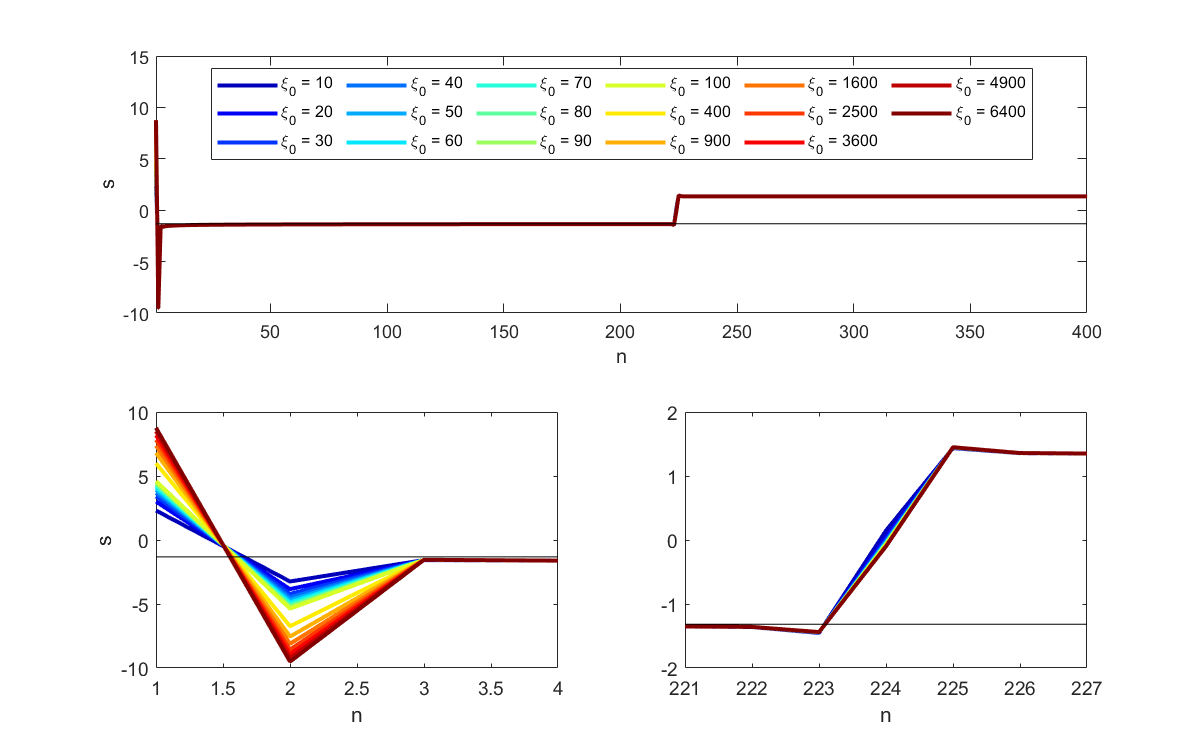}
\caption{Secant estimate of the slope of the log distance $\delta_n$ as a function of $n$, for Lew-Quarles orbits associated with different initial conditions $\xi_0$. Top panel: all curves are indistinguishable, suggesting that the convergence to the Freud orbit is uniform in $\xi_0$. The bottom row shows enlargements near $n = 1$ and $n = 224$, where differences between the curves are visible. In all panels, the black line corresponds to $s = \log(|\lambda_-|)$, where $\lambda_- = -2 + \sqrt 3$ is the stable eigenvalue of the linearization of dP1 near $P_\infty$.}
\label{fig:logdisslopes}
\end{figure} 

We witness the same turnaround behavior when replacing $\xi_0$ with different values ranging from 10 to 6,400. In Figure \ref{fig:logdisslopes}, we record the local secant approximation of the slope $s$ of the log distance $\delta_n$, for different values of $\xi_0$. In this case, all initial conditions were computed using 800 contractions. This figure strongly suggests that the rate of convergence of the Lew-Quarles orbits to the Freud orbit is uniform as all of the sub-graphs corresponding to different values of $\xi_0$ are effectively the same, giving the appearance of a single graph. Moreover, the thin black lines indicate that this rate is close to $|\lambda_-|$, where $\lambda_- = -2 + \sqrt 3$ is the stable eigenvalue of the linearization of dP1 about $P_\infty$. The bottom row shows regions where the slopes visibly depend on $\xi_0$. This occurs at the beginning of each orbit (small values of $n$, bottom left panel) and near the turnaround point (bottom right panel), when the accumulation of numerical errors starts to become noticeable and $s$ changes sign (for $n$ near 224 in the case of initial conditions calculated with 800 iterations of the Lew-Quarles contraction mapping, as is the case for Figure  \ref{fig:logdisslopes}).

\begin{figure}[ht]
\centering
\includegraphics[width=\textwidth]{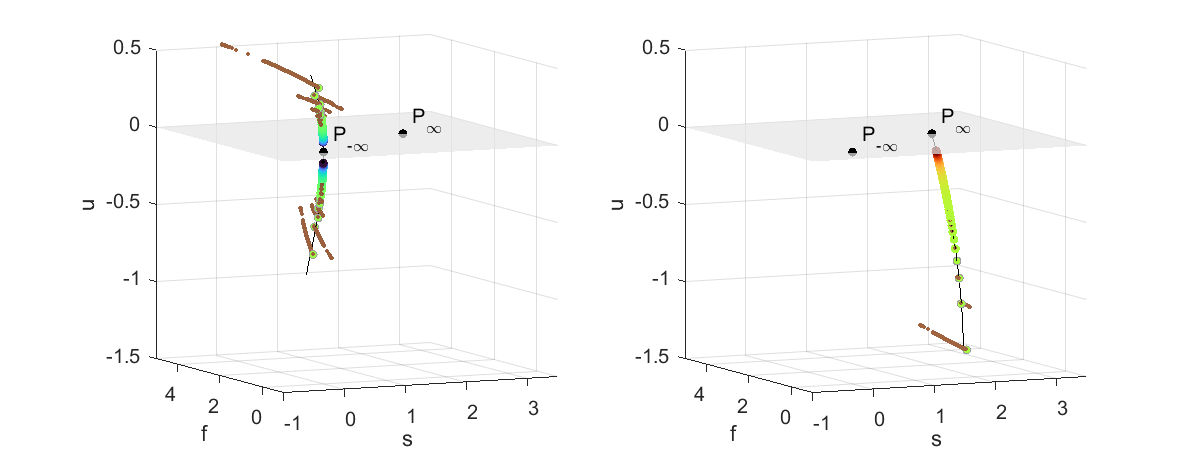}
\caption{Lew-Quarles orbits in $(s,f,u)$ coordinates, near $P_{-\infty}$ (left) and $P_\infty$ (right). Backwards iterates of $\mathcal S$ as well as the Freud orbit converge to $P_{-\infty}$ in a direction perpendicular to the plane $u=0$. Forward iterates quickly collapse onto the Freud orbit and converge to $P_\infty$ in a direction transverse, but not perpendicular, to the plane $u=0$. The solid curves represent approximations of the invariant curve given by equations \eqref{eq:CMP2} near $P_{-\infty}$ (left) and of the center manifold $F$ described by equations \eqref{eq:CMP1} near $P_\infty$ (right).}
\label{fig:sfuLQ}
\end{figure}

We now turn to a description of the Lew-Quarles solutions in the $(s,f,u)$ space. The right panel of Figure \ref{fig:sfuLQ} shows these orbits (in brown) as well as the Freud orbit (in color) in that space. This plot was obtained by applying the change of coordinates \eqref{eq:xyn_usf} to the initial conditions on $\mathcal S$ and their iterates, which were all numerically evaluated in the $(x,y,n)$ space.
In Appendix \ref{app:LQ_construction}, we extend Lew's and Quarles' results to arrive at the following theorem.
\begin{theorem} \label{th:LQ2}
All of the Lew-Quarles orbits, each defined by its initial condition on the set $\mathcal S$ of Theorem \ref{th:LQ}, converge to $P_\infty$ as $n \to \infty$.
\end{theorem}

\noindent It is natural to ask, at this stage, what dynamical systems theory can tell us about invariant manifolds containing $P_\infty$. In the autonomous limit of dP1 one knows that the phase space is foliated by invariant curves due to integrability. In the non-autonomous version such a foliation does not exist, but there are still invariant manifolds associated to fixed points such as $P_\infty$. We formulate the relevant results for us in this regard in the following two theorems which are consequences of the general center manifold theorem.
\begin{theorem} \cite{bib:i79} \label{th:CM}
Let $E_{cs}, E_u$ denote, respectively, the center-stable and unstable subspaces of $P_\infty$, and let $\Phi$ denote the 3-D map (\ref{eq:usf_sys}).
\begin{enumerate}
    \item Then there exists a smooth map $\chi: E_{cs} \to E_u$ whose graph $\mathcal{M}$ is tangent to $E_{cs}$ at $P_\infty$ and invariant under $\Phi$. Such an $\mathcal{M}$ is not necessarily unique, but we will refer to any such as a {\it center-stable} manifold for $\Phi$.
    \item Though $\mathcal{M}$ may not be unique, the coefficients of the Taylor series of $\chi$ are unique. 
    \item If $Q$ is a point in the phase space such that all its iterates $\Phi^n(Q)$ for $n$ larger than some $n_0 \in \mathbb{N}$ are in a certain fixed bounded neighborhood of $P_\infty$, then for any particular choice of center manifold $\overline{\mathcal M}$, $\text{dist}(\Phi^n(Q), \overline{\mathcal M}) \to 0$ as $n \to \infty$.
\end{enumerate}
\end{theorem}
\begin{theorem} \cite{bib:i79} \label{th:CM2}
Let $E_{c}, E_s$ denote, respectively, the center and stable tangent subspaces of $P_\infty$, in a local chart for a center-stable submanifold $\overline{\mathcal M}$ of $P_\infty$ and let $\widehat{\Phi}$ denote the restriction of $\Phi$ to $\overline{\mathcal M}$ in this chart.
\begin{enumerate}
    \item Then there exists a smooth map $\widehat{\chi}: E_{c} \to E_s$ whose graph $\mathcal{C}$ is a curve tangent to $E_{c}$ at $P_\infty$ and invariant under $\widehat{\Phi}$. Such a $\mathcal{C}$ is again not necessarily unique, but we will refer to any such as a {\it center} manifold (curve) for $\widehat{\Phi}$.
    \item Though $\mathcal{C}$ may not be unique, the coefficients of the Taylor series of $\widehat \chi$ are unique. 
    \item If $Q$ is a point on $\overline{\mathcal M}$ such that all its iterates $\widehat{\Phi}^n(Q)$ for $n$ larger than some $n_0 \in \mathbb{N}$ are in a certain fixed bounded neighborhood of $P_\infty$ in $\overline{\mathcal M}$, then $\text{dist}(\widehat{\Phi}^n(Q), \overline{\mathcal C}) \to 0$ as $n \to \infty$, for any particular choice $\overline{\mathcal C}$ of $\mathcal C$.
\end{enumerate}
\end{theorem}
The orbits described in Theorem \ref{th:LQ} satisfy the conditions of statement 3 of Theorem \ref{th:CM} and so are necessarily attracted to $\overline{\mathcal M}$. It is tempting to believe that in our case they are all actually {\it contained in} the same $\overline{\mathcal M}$ and determine it uniquely. Indeed it is difficult to imagine that a multiplicity of manifolds, comprised of orbits all limiting to
$P_\infty$ could coexist with structures of generic polar orbits. Unfortunately, the limitations of the center manifold theorem for maps do not enable us to conclude that on general grounds. (However, this is just one indication of the value in studying concrete examples of non-autonomous dynamics such as dP1.)

Despite these vagaries, the theoretical background provided by Theorems \ref{th:LQ}, \ref{th:LQ2}, and \ref{th:CM} together with the detailed numerical studies presented in this section will enable us to build a coherent and self-consistent picture of the dynamic state of affairs for non-polar orbits. In this picture ${\mathcal M}$ is unique and so contains all of the Lew-Quarles orbits. The latter orbits are then naturally interpreted as distinct center curves for $P_\infty$ as described in Theorem \ref{th:CM2}. In addition, the Freud orbit is chosen to represent the center manifold $\mathcal C$. 

We now look for an expansion of $\mathcal C$ in powers of $u$, by seeking a curve that is invariant under $\Phi$ and is tangent to the center direction of its linearization about $P_\infty$. We obtain
\begin{align}
s_\infty(u)&=2-u-\frac{\gamma}{6} u^2-\frac{\gamma}{36} u^3-\frac{\gamma  (3 \gamma +1)}{216} u^4+{\mathcal O}(u^{5}), \nonumber \\
& \label{eq:exp_sf_P1} \\
f_\infty(u)&=2-u+\frac{\gamma}{6} u^2+\frac{\gamma}{36} u^3-\frac{\gamma  (3 \gamma -1)}{216} u^4+{\mathcal O}(u^{5}) \nonumber
\end{align}
as $u \to 0.$ The iterative process leading to the above formulas can be continued to arbitrary order and expressions valid to 6\textsuperscript{th} order are provided in Equations \eqref{eq:CMP1} of Appendix \ref{app:CM}. These expressions are used to plot an approximation (black curve) of $F$ near $P_\infty$ in Figure \ref{fig:comb} and in the right panel of Figure \ref{fig:sfuLQ}. As shown in Figure \ref{fig:P1_expansions} of Appendix \ref{app:CM}, they capture the Freud orbit extremely well, even for values of $u$ of order one, thereby providing numerical support to our conjecture that the Freud orbit $F$ is well approximated by the center manifold $\mathcal C$ of $P_\infty$. We note that by the center manifold theorem the above expansions are independent of the particular choice of $\mathcal{C}$.

\subsection{Backward iterates}
\label{sec:backward_iterates}
The map sequence $\phi_n$ given by (\ref{eq:napdintrotwo}) has a sequence of inverse maps $\psi_n$ given by 
\begin{align} \label{eq:INV} 
x_n&=   y_{n+1}\\ 
y_n&= \dfrac{n}{N r y_{n+1}}-\dfrac{1}{r}-x_{n+1}-y_{n+1}.    \nonumber
\end{align} 
 The map $\psi_n$ is a birational mapping singular along the $x$-axis. We are interested in studying extensions of the above non-polar orbits in reverse time using $\psi_n$.  In this case the singularities enter our consideration since the initial point of the Freud orbit, $(x_1, y_1 = 0)$ corresponding to $\xi_0 = 0$,  lies on the $x$-axis. This is the only Lew-Quarles orbit that has this issue. Both the $\phi_n$ and $\psi_n$ mappings can be extended to be defined along their respective singular axes by a resolution of singularities process (see \cite{bib:kny17}). However we take a simpler approach that is more relevant to our asymptotic change of variables

We address the singularity of the map \eqref{eq:INV} on the $y=0$ axis  by taking a limit from the backward iterates of the Lew-Quarles orbits for which $x_0 =\xi_0 > 0$. In that case, $y_1>0$ and the term $n / (N r y_1)$ in the inverse mapping $\psi_0$ is well defined and equal to $0$ since $n=0$. Therefore, for $y_{1}=0$ and $n=0$, we define $y_0$ in the mapping $\psi_0$ by
\[
y_0 = 0 -\dfrac{1}{r}-x_{1}-y_{1} = -\dfrac{1}{r}-x_{1}.
\]

Going back to Figure \ref{fig:xyLQ}, $\psi_0({\mathcal S})$, which is well defined since $y_1 > 0$, is the collection of points shown in brown in the fourth quadrant. As mentioned above, this set of points limits to a point on the $y$-axis, which we defined to be the image of $(x_1^F,y_1^F)$ under $\psi_0$. Successive applications of $\psi_n$, $n < 0$, to $\psi_0(x_1^F,y_1^F)$ define the backward iterates $\left\{\psi_n \circ \psi_{n+1} \circ \cdots \circ \psi_0(x_1^F,y_1^F)\right\}$ of the Freud orbit. The Freud orbit is shown in green and blue as $n$ becomes more negative, while backward iterates of $\mathcal S$ under $\psi_n, n \le -1$ are in brown. They alternate between the second and fourth quadrants. As was the case for its forward iterates, the backward iterates of $\mathcal S$ converge to the Freud orbit as $n$ becomes more negative. The numerically estimated rate of convergence appears to be sub-exponential, which is consistent with all of the eigenvalues of the linearization of dP1 at $P_{-\infty}$ being on the unit circle. Backward iterates of $\mathcal S$ in the fourth quadrant are not visible because they are superimposed with $\psi_0({\mathcal S})$ when projected onto the $(x,y)$ plane. Because the dP1 system \eqref{eq:napdintrotwo} is singular along the Freud orbit when $n=0$ (since $x_0 = 0$) but not along the other Lew-Quarles orbits, the above construction allows us to view the former as a singular limit of the latter.

As illustrated in Figure \ref{fig:sfuLQ}, our numerical simulations show that in the $(s,f,u)$ coordinates, backward iterates of $\mathcal S$ converge to the fixed point $P_{-\infty}$ in a direction perpendicular to the invariant plane $u=0$. Recall that the linearization of \eqref{eq:usf_sys} about $P_{-\infty}$ has eigenvalues $-1$ and $\pm i$, and the eigendirection associated to $-1$ is perpendicular to the plane $u=0$. To better understand the dynamics near $P_{-\infty}$, we look for an invariant curve parametrized by $u$, of the form
\[
s=s_{-\infty}(u) = \sum_{k=2}^\infty a_k u^k, \qquad f=f_{-\infty}(u) = \sum_{k=2}^\infty b_k u^k.
\]
Requiring that $s_{n+1} = s_{-\infty}(u_{n+1})$ and $f_{n+1} = f_{-\infty}(u_{n+1})$ be satisfied when ($s_n = s_{-\infty}(u_n)$, $f_n=_{-\infty}(u_n)$, $u_n$) and $(s_{n+1},f_{n+1},u_{n+1})$ are related by the mapping \eqref{eq:usf_sys}, leads to a consistency relation in powers of $u$, which in turn defines a set of equations for the coefficients $a_k$ and $b_k$. This procedure gives
the following formal expressions for $s_{-\infty}$ and $f_{-\infty}$ as $u \to 0$ near $P_{-\infty}$,
\begin{align}
s=& s_{-\infty}(u) = -\frac{\gamma}{2} u^2 - \frac{\gamma}{4} u^3 + \frac{\gamma}{8} (\gamma - 1) u^4 + {\mathcal O}(u^{5}), \nonumber \\
& \qquad \label{eq:exp_sf} \\
f=& f_{-\infty}(u) = +\frac{\gamma}{2} u^2 + \frac{\gamma}{4} u^3 + \frac{\gamma}{8} (\gamma + 1) u^4 + {\mathcal O}(u^{5}).\nonumber
\end{align}
The above relationships can be pushed to higher order with the help of computer algebra software. Equations \eqref{eq:CMP2} in the appendix, valid to order 10, are used to plot the thin black invariant curve through $P_{-\infty}$ displayed in Figures \ref{fig:comb} and \ref{fig:sfuLQ}. They capture the Freud orbit very well near $P_{-\infty}$, even for relatively large values of $u$, as further illustrated in Figure \ref{fig:P2_expansions} of Appendix \ref{app:CM}.

In summary, the combination of numerical and analytical investigations presented in this section suggests the following picture: in $(s,f,u)$ coordinates, the Lew-Quarles orbits form a family of heteroclinic connections between $P_{-\infty}$ and $P_\infty$. All iterates are defined for all values of $n$, but $u_0 \to -\infty$ (equivalently $x_0 \to 0$) as $\xi_0 \to 0$, where $\xi_0$ parametrizes the family of Lew-Quarles orbits. This limit defines the Freud orbit, which leaves $P_{-\infty}$ along an invariant curve perpendicular to the plane $u=0$ and converges to $P_\infty$ along its center direction.

\section{Characterization of the Freud orbit}
\label{sec:Freud}

This section provides the conceptual framework for understanding a striking observation: invariant curves near the fixed points of \eqref{eq:usf_sys} {\it track} periodic points of the associated autonomous limits of dP1. This conceptual framework is grounded in a novel process that realizes a simple and elegant mechanism for generating classical asymptotic expansions of the Freud orbit. 

In Section \ref{sec:NP}, we made the conjecture that near $P_{-\infty}$ and $P_\infty$, the Lew-Quarles orbits lie on invariant curves whose expansions are given by \eqref{eq:CMP2} and \eqref{eq:CMP1} respectively. We now use this assumption to determine the behavior of $u_n$ as a function of $n$ along the Freud orbit near $P_{-\infty}$ and $P_\infty$.
Near $P_{-\infty}$, we look for a sequence $\{u_n\}$ that converges to $0$ as $n \to -\infty$ and is consistent with the last equation of \eqref{eq:usf_xyn}. We thus require that
\begin{equation}
\frac{n}{N} = \frac{s_{- \infty}(u_n) + f_{- \infty}(u_n) + u_n - 1}{r u_n^2} \Longleftrightarrow \gamma n u_n^2 - u_n + 1 = s_{- \infty}(u_n) + f_{- \infty}(u_n).
\label{eq:P2_un}
\end{equation}
Since per \eqref{eq:exp_sf} $s_{- \infty}(u_n) + f_{- \infty}(u_n) = {\mathcal O}(u_n^{4})$, we see that at dominant order, the iterates $\{u_n\}$ solve $\gamma\, n\, u_n^2 - u_n + 1 = 0$, which is the equation defining the period-2 solutions of the autonomous dP1 system in $(s,f,u)$ coordinates. As further described in Appendix \ref{app:Aut}, this system is obtained from dP1 by setting $\alpha = n/N$ and then assuming that $\alpha$ is a constant parameter. We call the resulting autonomous map $\alpha$-dP1. Figure \ref{fig:p2pts} of Appendix \ref{sec:period2points} illustrates the numerical convergence of the Freud orbit to the sequence of period-two points of $\alpha$-dP1 defined in Equation \eqref{eq:peralpha}, as $n \to -\infty$. Equation \eqref{eq:P2_un} also provides an expansion of $u_n$ as a function of $n$,
\begin{equation}
\label{eq:unP2}
u_{-\infty,n} = u_{-\infty,n}^\pm = \pm \frac{1}{\sqrt{- \gamma n}} + \frac{1}{2 \gamma n} \pm \frac{1}{8 (- \gamma n)^{3/2}} + {\mathcal O}\left((-n)^{-5/2}\right) \quad \text{as } n \to -\infty,
\end{equation}
which is obtained by writing $u_{-\infty,n}$ as a Laurent expansion in powers of $\sqrt{-n}$, substituting into the right-hand equation of \eqref{eq:P2_un}, and solving term by term. At this point, this expansion is formal because we do not have a proof of the existence of the smooth functions $s_{-\infty}$ and $f_{-\infty}$ that appear in expansions \eqref{eq:exp_sf}, and therefore do not have enough control to bound the remainder in \eqref{eq:unP2}. The approach however, is very general. The definition of $\alpha$ in terms of $s$, $f$, and $u$, together with Equation \eqref{eq:CMP2}, immediately leads to two results: that solutions of dP1 on the invariant curve associated with \eqref{eq:CMP2} track the sequence of period-two points of the autonomous dP1 system as $n \to -\infty$, and the expansion \eqref{eq:unP2}. 

\begin{figure}[ht]
\centering
\includegraphics[width=.9\textwidth]{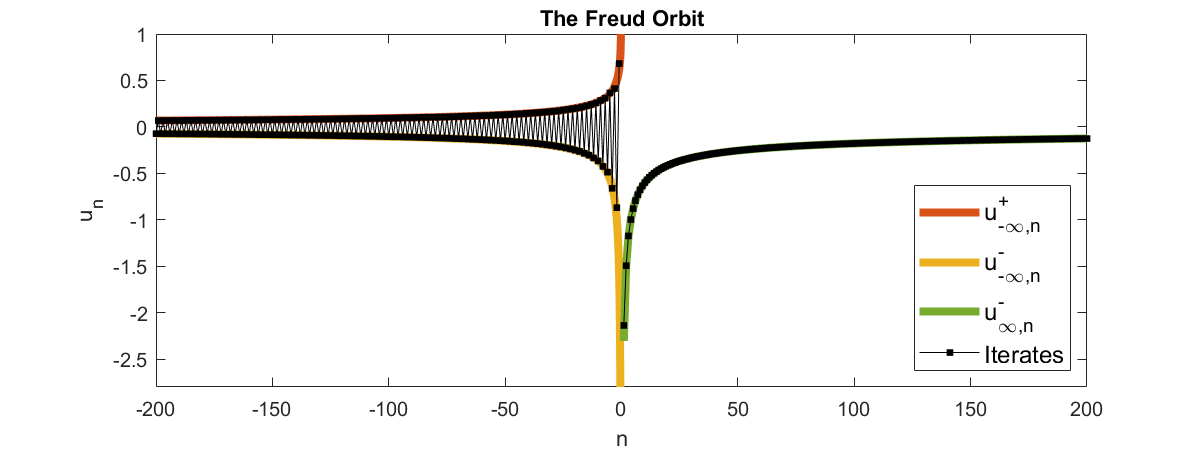}
\caption{Values of $u_n$ for iterates of the Freud orbit (black, connected dots), together with the asymptotic expansions $u_{-\infty,n}^\pm$ defined in \eqref{eq:unP2} for $n < 0$ (solid red and yellow curves), and $u_{\infty,n}^-$ defined in \eqref{eq:unP1} for $n > 0$ (solid green curve).
\label{fig:Comb3}}
\end{figure}

Near $P_\infty$, we proceed in a similar fashion, although in that case the existence of the center manifold $\mathcal C$ makes the resulting expansions asymptotic. Equation \eqref{eq:P2_un} is replaced by
\[
s_{\infty}(u_n) + f_{\infty}(u_n) = \gamma\, n\, u_n^2-u_n+1.
\]
Using \eqref{eq:exp_sf_P1}, we see that at leading order, the iterates $\{u_n\}$ solve $\gamma\, n\, u_n^2+u_n-3= 0$, which is the equation defining the fixed point solutions of the autonomous dP1 system in $(s,f,u)$ coordinates (see Appendix \ref{app:Aut}). Solving for $u_n$ as a function of $n$, we find

\begin{equation}
\label{eq:unP1}
u_{\infty,n} = u_{\infty,n}^\pm = \pm \sqrt{\frac{3}{\gamma n}}-\frac{1}{2 \gamma  n} \pm \frac{1}{\left(8 \sqrt{3}\right) (\gamma n) ^{3/2}} + {\mathcal O}\left(n^{-5/2}\right), \quad \text{as } n \to \infty.
\end{equation}
We note that control of the big $\mathcal O$ term follows from the center manifold theorem that guarantees the existence of all Taylor coefficients of $s_\infty$ and $f_\infty$, and then a straightforward application of Taylor's remainder theorem. Plots of $u_n$ as a function of $n$ for the numerically estimated Freud orbit, and of the expansions $u_{-\infty,n}^\pm$ and $u_{\infty,n}^-$ given above are provided in Figure \ref{fig:Comb3}; they show excellent agreement, even for values of $n$ near 0. This figure dramatically reinforces the fact that the Freud orbit is singled out among the Lew-Quarles orbits by its singularity at $n=0$.

Given their uniform exponential convergence to the Freud orbit exhibited above, it is natural to further conjecture that 
all the Lew-Quarles orbits defined in Theorem \ref{th:LQ} have this same asymptotic expansion, and so differ only in terms that are beyond all orders with respect to the algebraic gauge $n^{-1/2}$.

\section{Conclusions}
\label{sec:conclusion}

By introducing the asymptotic change of coordinates \eqref{eq:xyn_usf}, we transformed the dP1 mapping \eqref{eq:napdintrotwo} into the 3-dimensional discrete autonomous dynamical system \eqref{eq:usf_sys}, whose two fixed points, $P_{-\infty}$ and $P_\infty$, correspond to solutions $\{x_n\}$ of dP1 that grow without bounds as $|n| \to \infty$. This transformation, together with a combination of analytical and numerical investigations, allowed us to characterize known non-polar orbits of dP1 as heteroclinic connections between these two fixed points. In particular, we described the Freud orbit as a singular limit of the Lew-Quarles orbits. By understanding how these solutions leave $P_{-\infty}$ to converge to $P_\infty$ as $n$ increases, we discovered that they track sequences of points constructed from period-2 (near $P_{-\infty}$) and period-1 (near $P_\infty$) points of the autonomous counterpart of dP1. Moreover, our results are consistent with and in many aspects support the conjecture that the Lew-Quarles orbits are on center-stable manifolds of $P_\infty$ and, as $n \to \infty$, exponentially converge to the Freud orbit, which itself lives on a center manifold of $P_\infty$. The presence of invariant curves that contain the Lew-Quarles (as $n \to \infty$) and Freud orbits provides a method to find explicit expansions of these solutions in powers of $|n|^{-1/2}$ as $|n| \to \infty$. These expansions are asymptotic near $P_\infty$ and formal near $P_{-\infty}$.

Compared to their continuous counterparts, many aspects of discrete dynamical systems remain relatively unexplored, especially in the non-autonomous case. In this environment, the study of concrete examples takes on an added value. Due to its analytically completely integrable discrete limit, dP1 is of particular interest. The resulting rich inherent structure therefore provides an auspicious environment for the exploration of explicit features of non-autonomous discrete dynamics that cannot be directly attacked in a general setting. This paper does that in the context of what we have termed non-polar orbits and relates its findings to previously known aspects of Painlev\'e systems. This approach also leads to interesting applications regarding the behavior of specific orbits as $|n| \to \infty$.

The autonomous dP1 system has orbits that are degenerations of its  quasi-periodic elliptic solutions and correspond to trigonometric solutions along a separatrix. This separatrix, and the solutions along it, were examined in detail from a dynamical systems perspective in \cite{bib:be19}  and used to better understand their relation to the combinatorial problem of {\it geodesic distance} in the  enumeration of planar graphs. That paper also examined a related system, $x_n (x_{n+1} + x_{n-1} + 1/r) =n / Nr$, which has relevance to the enumeration of labelled trees and super-Brownian excursions. Both of these are instances of the application of dP1 type systems to other areas of mathematics. The present manuscript reveals a subtle connection between the separatrices of autonomous dP1 discussed in \cite{bib:be19} and Freud's orbit. Specifically, we discover a novel and deep connection between the mapping \eqref{eq:napdintrotwo} and its autonomous counterpart, as the sequence $\{x_n\}$ in Freud's recurrence follows the fixed points of a collection of autonomous mappings. In Appendix \ref{app:Aut}, we provide a geometric description of the important role played by these separatrices in making it possible for the Freud orbit to grow without bounds as a solution of dP1. In addition, the principal message of this article is that Freud's orbit coincides with a center manifold of a fixed point at infinity. This property of having, effectively, a limiting fixed point at infinity is highly atypical. We do not believe that such a connection between autonomous fixed points and an asymptotic non-autonomous fixed point has been noticed before in the literature. 
  
Another novel discovery is the realization that recursive constructions of invariant curves provide an elegant mechanism for generating asymptotic expansions of non-polar solutions of dP1 as $|n| \to \infty$ (see  \eqref{eq:unP2} and \eqref{eq:unP1}). For $n > 0$, these have relevance for another combinatorial problem related  to random tilings of Riemann surfaces (or geometric foams and quantum gravity in the physics literature). Together with the results of \cite{bib:be19}, we therefore now have three examples of dynamical systems structures (specifically two instances of a stable manifold of a hyperbolic fixed point in \cite{bib:be19}, and a center manifold here) associated with solutions to combinatorial problems. Such connections between different areas of mathematics are highly intriguing. In a subsequent paper we will develop some of these connections using our recursive construction of invariant curves and their associated asymptotic expansions. Methods developed in \cite{bib:tip20} will enable us to demonstrate the key feature of uniform validity in $r$ and  $N$, within appropriate ranges, for these expansions.

As mentioned in the introduction, M\'{a}t\'{e}, Nevai, and Zaslavsky in 1985 \cite{bib:mnz85}  established the  existence of an asymptotic expansion for the Freud orbit in powers of $n^{1/2}$, whose leading behavior can be shown to be \cite{bib:tip20}
\begin{align} \label{eq:introlikewn}
x_n& = \sqrt{\dfrac{n}{3rN}}-\dfrac{1}{6r}+\dfrac{\sqrt{12}}{144}\sqrt{\dfrac{N}{n}}\dfrac{1}{r^{3/2}} + O(1/n).
\end{align}
We also note that in \cite{bib:jos97} Joshi carried out a formal Painlev\'e dominant balance analysis, seeking asymptotic fixed points of dPI in terms of an asymptotic constant of motion. From the perspective of this paper one can see that this formal approach correctly captured the leading order behavior at $P_\infty$ and $P_{-\infty}$. It would be interesting to determine if and where this approach might match the ``level set" structure evident in Figure \ref{fig:u=0}, and possibly even approximate our exact trajectories in the vicinity of $P_{\pm\infty}$.
The mechanisms just referred to in the previous paragraph reproduce and extend this type of asymptotic analysis. Indeed, in these connections, our results may say more about the asymptotics of orthogonal polynomials. 

Our subsequent work will also make central use of the realization that the leading order term in the large $n$ asymptotic expansion of \cite{bib:emp08} in fact coincides with our expression for the fixed points of autonomous dP1. This, together with the analysis developed in Section \ref{sec:Freud},
enables us to devise an efficient scheme for explicitly counting topologically
distinct quadrangulations, involving a fixed number of tiles, of compact Riemann surfaces. Other, higher dimensional, Painlev\'e systems have orbits of Freud type corresponding to degree $2 \nu$ potentials that in turn yield information about the enumeration of more general polygonal tilings. For potentials of odd dominant degree traditional methods of orthogonal polynomial theory break down and one must consider generalizations such as non-Hermitian orthogonal polynomials \cite{bib:ep12}, which presents obstacles to asymptotic analysis. The particular case of a cubic potential, which relates to topological triangulations of surfaces, is of special interest and involves a non-autonomous planar Painlev\'e dynamical system. The dynamical systems approach developed in this paper may help to overcome some of these obstacles.

Finally, we note that our work has focused on the dP1 regime with $r >0$. However, there is also interest in the singular regime where $r \to -1/12$, the so-called ``Boutroux regime'' that has been explored in \cite{bib:jl15}  and \cite{bib:dk06}.  This is the so-called double scaling limit of random matrix theory, which provides a bridge between the discrete and continuous versions of Painlev\'e I. A problem of general physical importance is to study and relate the behavior of a correlation function for a statistical mechanical system at small values of a scaling parameter  to its behavior at large values. This is called a {\it connection problem}. For a class of explictly solvable statistical mechanical problems this is related to connection problems for continuous Painlev\'e equations. For Painlev\'e I and II a systematic study of this problem was carried out in \cite{bib:jk}.
A natural extension of our work is to explore the potentially interesting connection problem between the limits $r \to \infty$ and $r \to -1/12$ from a dynamical systems perspective.

As stated at the outset of these Conclusions, the work presented here combines theory and asymptotics with high-precision numerical simulations to arrive at a detailed and compelling picture for the structure of non-polar dynamics in the dP1 system. The main work needed to convert this picture into fully rigorous statements centers on establishing our stated conjectures about center-stable and center manifolds laid out in Section \ref{sec:NP}. This challenge is of independent interest. Its resolution will have significance for discrete dynamical systems theory broadly speaking. Within that scope particular open questions of importance include
\begin{enumerate}
    \item The rigorous verification of the numerically evident fact that all of the Lew-Quarles orbits converge exponentially to the Freud orbit at a uniform rate as $n \to \infty$.
    \item A more general conceptual understanding of the unique singularity formation in the
    Freud orbit as compared to the other Lew-Quarles orbits. In other (continuum) examples of center manifolds (see e.g. \cite{bib:meiss07} Section 5.6) one observes that singularities form in finite time along those center manifolds whose asymptotic expansions do not have any corrections beyond all orders. Could this be the case for Freud's orbit?
\end{enumerate}

\appendix

\section{Orthogonal polynomials}
\label{app:OP}
We refer the reader to \cite{bib:tip20} for a thorough treatment of the orthogonal polynomials given by the weight \eqref{eq:orthoweight}. The following is a condensed derivation for an exact integral representation of the Freud initial condition. We denote the weight, recurrence, and moments by
\begin{equation}
w(\lambda)=\exp[-N(\lambda^2/2+(r/4)\lambda^4)]
\label{eq:orthoweight}
\end{equation}
\begin{equation}
\lambda p_n(\lambda) = b_{n+1}p_{n+1}(\lambda)+b_np_{n-1}(\lambda),
\label{eq:generalrec}
\end{equation} 
\begin{equation*}
\mu_i=\int_{-\infty}^\infty \lambda^i w(\lambda) d\lambda.
\end{equation*} 
We look for orthonormal polynomials, so computing $p_0$ is straightforward:
\[
1= \int_{-\infty}^\infty p_0^2 w(\lambda) d\lambda = p_0^2(\lambda) \int_{-\infty}^\infty w(\lambda) d\lambda = p_0^2 \mu_0 \Longrightarrow
p_0=\dfrac{1}{\mu_0^{1/2}}.
\]
To find $p_1$, we apply the recurrence \eqref{eq:generalrec}: 
$$\dfrac{ \lambda}{\mu_0^{1/2}} =\lambda p_0 = b_1 p_1(\lambda) \Longrightarrow p_1(\lambda)=\dfrac{\lambda}{b_1\mu_0^{1/2}}$$
and normalize $p_1$ to solve for $b_1$:
$$\frac{1}{b_1^2\mu_0}\int_{-\infty}^\infty \lambda^2 w(\lambda) d\lambda =1 \Longleftrightarrow \frac{\mu_2}{b_1^2\mu_0}=1.$$
Thus,
\[
b_1^2 = \frac{\mu_2}{\mu_0}, \qquad
p_1(\lambda) =\dfrac{\lambda}{\mu_2^{1/2}}.
\]
For $p_2$, again applying the recurrence \eqref{eq:generalrec} leads to
$$\dfrac{ \lambda^2}{\mu_2^{1/2}} =\lambda p_1(\lambda) = b_2 p_2(\lambda) +b_1 p_0 \Longrightarrow p_2(\lambda) = \dfrac{1}{b_2}\left(\dfrac{\lambda^2}{\mu_2^{1/2}}-\dfrac{\mu_2^{1/2}}{\mu_0} \right).$$
Then, normalizing $p_2$ leads to
$$\dfrac{1}{b_2^2}\int_{-\infty}^\infty \left(\frac{\lambda^4}{\mu_2}-\frac{2 \lambda^2}{\mu_0}+\frac{\mu_2}{\mu_0^2}\right)w(\lambda)d\lambda=1 \Longleftrightarrow \dfrac{1}{b_2^2} \left( \frac{\mu_4}{\mu_2}-\frac{2\mu_2}{\mu_0}+\frac{\mu_2}{\mu_0}\right)=1.$$
Solving for $b_2^2$, we find
\begin{equation*}
b_2^2=\dfrac{\mu_4\mu_0-\mu_2^2}{\mu_0\mu_2}, \qquad
p_2(\lambda)=\lambda^2\left(\dfrac{\mu_0\mu_2^{1/2}}{\mu_4\mu_0-\mu_2^2}\right)-\dfrac{\mu_2^{3/2}}{\mu_4\mu_0-\mu_2^2}.
\end{equation*} 
Summarizing the above, we have
$$b_1^2=\dfrac{\mu_2}{\mu_0},\ \ \ b_2^2=\dfrac{\mu_4\mu_0-\mu_2^2}{\mu_0\mu_2},$$
which defines the point $(x_2=b_2^2,x_1=b_1^2)$ on the Freud orbit in terms of moments of the weight function $w$.

We note that we can also define $b_0 = 0$, thereby obtaining an equivalent way to initialize the $b_n$ sequence. Indeed, using Freud's equation \eqref{eq:Freudeq}, we have
\[
N(rb_1^2(b_{2}^2+b_{1}^2+b_{0}^2)+b_1^2) = 1 = N\left(r\dfrac{\mu_2}{\mu_0}\left(\dfrac{\mu_4\mu_0-\mu_2^2}{\mu_0\mu_2}+\dfrac{\mu_2}{\mu_0}\right)+\dfrac{\mu_2}{\mu_0}\right).
\]
Verification of the second equality is a simple application of integration by parts. 

\section{Lew-Quarles construction}
\label{app:LQ_construction}
The proof of Theorem \ref{th:LQ} is based on a contraction argument. What we present here is a mild modification, applicable for system \eqref{eq:napdintrotwo}, of what is developed in \cite{bib:lq83}. A different proof that does not involve a contraction mapping argument is provided in Corollary 5.7 of \cite{bib:ansva15}.

Let $\{\xi_n\}$ be a solution of dP1 that remains in the first quadrant. In terms of \eqref{eq:napdintrotwo}, we have $x_n = \xi_n$ for $n > 0$, and $y_1 = \xi_0$. It is straightforward to rewrite \eqref{eq:napdintrotwo} as
\bea \label{LQquad} 
\left(\frac{\xi_n}{\kappa_n}\right)^2 &=& 1 -  \frac{\xi_n}{\kappa_n} \left(\frac{\xi_{n-1} + \xi_{n+1}}{\kappa_n}  + \frac1{r \kappa_n}\right),
\eea
where $\kappa^2_n = \frac{n}{Nr}$.  Setting 
\bea
\sigma &=& \frac{\xi_n}{\kappa_n}\\
\tau &=& \frac{\xi_{n-1} + \xi_{n+1}}{2 \kappa_n}  + \frac1{2 r \kappa_n},
\eea
we see that (\ref{LQquad}) has the form of a simple quadratic equation in $\sigma$, $\sigma^2 + 2 \tau \sigma -1 = 0$. Since we have assumed $\xi_n$
(and $\kappa_n$ as well) to be positive, $\sigma$ may be uniquely expressed as the positive root of the quadratic
\bea
g(\tau) &=& -\tau + \sqrt{1 + \tau^2} > 0.
\eea
In terms of this $g$ one defines a mapping, $T$ on the space of non-negative sequences $x = (\xi_1, \xi_2, \dots)$ for fixed 
$\xi_0 \geq 0$ and $c = (\kappa_1, \kappa_2, \dots)$ given by
\bea \label{contraction}
(T x)_n &=& \kappa_n \cdot g\left(  \frac{\xi_{n-1} + \xi_{n+1}}{2 \kappa_n}  + \frac1{2 r \kappa_n} \right), \,\,\,\,\, n \geq 1.
\eea
This mapping is a contraction which converges, component-wise, to a unique, positive fixed point which, by the above construction is a solution of  (\ref{eq:napdintrotwo}). This contraction also provides a highly efficient numerical algorithm for calculating non-polar solutions of dP1.

Lew and Quarles \cite{bib:lq83} also provide the essential ingredients to prove Theorem \ref{th:LQ2}. The first of these ingredients is a basic estimate.
\begin{lemma}
If $Tx = x$, then 
$$Tc < x <c.$$
\end{lemma}
Using this lemma as well as the structure of the map \eqref{eq:napdintrotwo} one can show that
\begin{proposition}
If $Tx=x$, then $\dfrac{\xi_m}{\kappa_m}$ has a limit as $m \rightarrow \infty$, and 
$$\theta = \lim_{m \to \infty} \frac{\xi_m}{\kappa_m} = (\mu+1+ \mu^{-1})^{-1/2},$$
where $\displaystyle \mu = \lim_{m \to \infty} \frac{\kappa_{m+1}}{\kappa_m}$, which is equal to 1.
\end{proposition}
It is then straightforward to see that $\xi_n \sim \sqrt{\frac{n}{3Nr}}$. Theorem \ref{th:LQ2} follows, since
\begin{align*}
s_m & = \frac{\xi_{m-1}}{\xi_m} + 1 + \frac{1}{r \xi_m} = \frac{\xi_{m-1}}{\kappa_{m-1}}\frac{\kappa_{m-1}}{\kappa_m}\frac{\kappa_m}{\xi_m} + 1 + \frac{1}{r \kappa_m} \frac{\kappa_m}{\xi_m} \to 1 + 1 + 0 = 2,\\
f_m & = \frac{\alpha}{r \xi_m^2} - \frac{\xi_{m-1}}{\xi_m} = \frac{m}{N r \kappa_m^2}\frac{\kappa_m^2}{\xi_m^2} - \frac{\xi_{m-1}}{\xi_m} \to 1 \cdot 3 - 1 = 2,\\ 
u_m &= - \frac{1}{r \xi_m} \to 0
\end{align*}
as $m \to \infty$. We note that $\xi_n \sim \sqrt{\frac{n}{3Nr}}$ matches the asymptotic expansion \eqref{eq:introlikewn} for Freud's orbit.
\medskip

From a broader perspective (\ref{contraction}) may be viewed as a family of infinite dimensional mappings parametrized by the sequences $(\kappa_1, \kappa_2, \dots)$. To demonstrate that this mapping has a fixed point that is unique in the space of non-negative sequences $x = (\xi_1, \xi_2, \dots)$, all that the Lew-Quarles proof requires is that
\begin{enumerate}
    \item $\xi_0 \geq 0$;
    \item the $\kappa_n$ are all positive;
    \item $\inf_m \kappa_m/m =0$.
\end{enumerate}
These conditions clearly hold for the original Lew-Quarles sequences. But now suppose one has a solution of dP1 such that there is an $n_0$ where $x_n > 0$ for all $n \geq n_0$. Then set $\xi_0 = x_{n_0}$ which is positive and so satisfies the first condition above. Define a new sequence $(\kappa_{n_0 + 1}, \kappa_{n_0 + 2}, \dots)$ where the $\kappa_j$ are those defined earlier. This {\it truncated} defining sequence clearly satisfies the remaining two conditions above since  $\inf_m \kappa_{m+ n_0}/m =0$. So this new mapping must have a unique fixed point by the original Lew-Quarles argument. But such a fixed point is by construction a solution to dP1 with initial value passing through $x_{n_0}$ and remaining positive thereafter. By the analogue of Theorem \ref{th:LQ2}, this additional orbit limits to $P_\infty$.

\section{Fixed Points and Period 2 Points of Autonomous dP1}
\label{app:Aut}

\subsection{$\alpha$-dP1 Hyperbolic Fixed Point for $\alpha > 0$}
\label{sec:hyperbolicfixedpoint}

\begin{figure}[bht]
\centering
\includegraphics[width=\textwidth]{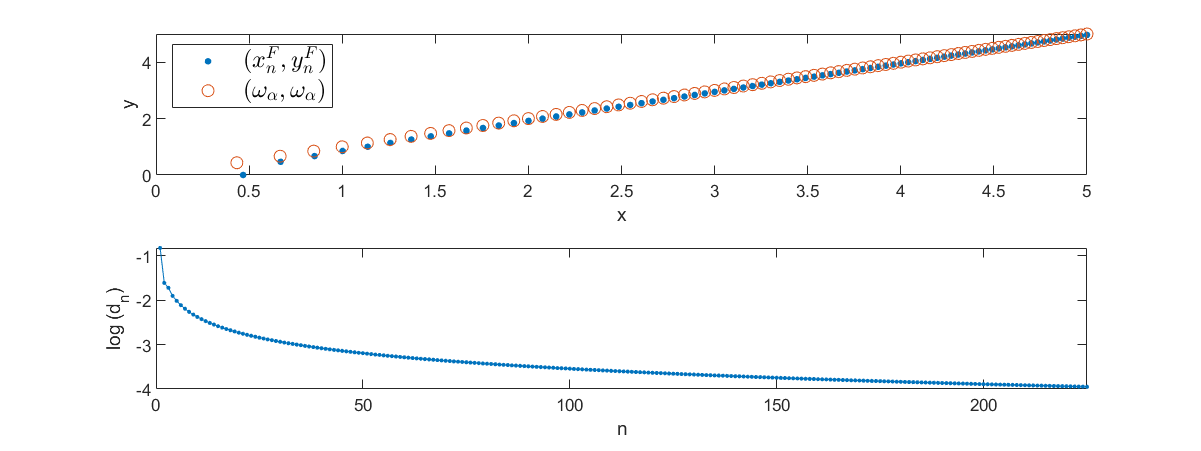}
\caption{Top: Points on the Freud orbit (dots) and in the sequence $W_n$, for $1 \le n \le 80$ and $N = 1$. Bottom: Log-distance $\log\left(d_n \right)$ as a function of $\alpha = n$.}
\label{fig:hypfixedpoints}
\end{figure}

The mapping introduced in Equation \eqref{eq:napdintrotwo} is based on the family of what we call $\alpha$-dP1 mappings:
\begin{align}
\overline{x}&= \dfrac{\alpha}{rx}-\dfrac{1}{r}-x-y \nonumber \\
\overline{y}&= x \label{eq:alpha},
\end{align}
in which $\alpha$ is the extension of $n/N$ to a continuous parameter. These mappings have two fixed points for $\alpha > 0$, one hyperbolic and the other elliptic.  In $(x, y, \alpha)$ coordinates the hyperbolic points are given by
\begin{equation}
(\omega_{\alpha},\omega_{\alpha}) :=\left(\dfrac{-1+ \sqrt{1+12\alpha r}}{6r},\dfrac{-1+ \sqrt{1+12\alpha r}}{6r}\right).
\label{eq:alphafixed}
\end{equation}
In $(s, f, u)$ coordinates, they are
\begin{equation}
(s_\alpha, f_\alpha, u_\alpha) :=\left(2 - u_\alpha, 2 - u_\alpha, -6\left(-1 + \sqrt{1 +12 r \alpha}\right)^{-1}\right),
\label{eq:alphafixed2}
\end{equation}
with $u_\alpha$ solving the quadratic equation $\alpha\, r\, u_\alpha^2 + u_\alpha - 3 = 0$. The top panel of Figure \ref{fig:hypfixedpoints} shows part of the sequence of points $W_n = \left\{(\omega_{\alpha},\omega_{\alpha}), \alpha = n/N \right\}_{n \in \mathbb{N}}$ (open circles). In \cite{bib:tip20} it was first noted that this plot looks remarkably like the plot of the Freud orbit (dots in the top panel of Figure \ref{fig:hypfixedpoints}) and then analytically established that, in fact, the sequence $W_n$ approximates the Freud orbit through order 2 in powers of $n^{-1/2}$. The bottom panel of Figure \ref{fig:hypfixedpoints} shows a plot of the log-distance $\log(d_n)$, where
\[
d_n = \sqrt{(x_n^F - \omega_\alpha)^2 + (y_n^F - \omega_\alpha)^2}, \quad \alpha = \frac{n}{N}
\]
for $N = 1$ and $n \in \{1, \cdots, 225\}$, illustrating the convergence of the two sequences of points to one another as $n \to \infty$.

\subsection{$\alpha$-dP1 Period 2 Points for $\alpha < 0$}
\label{sec:period2points}

\begin{figure}
\centering
\includegraphics[width=\textwidth]{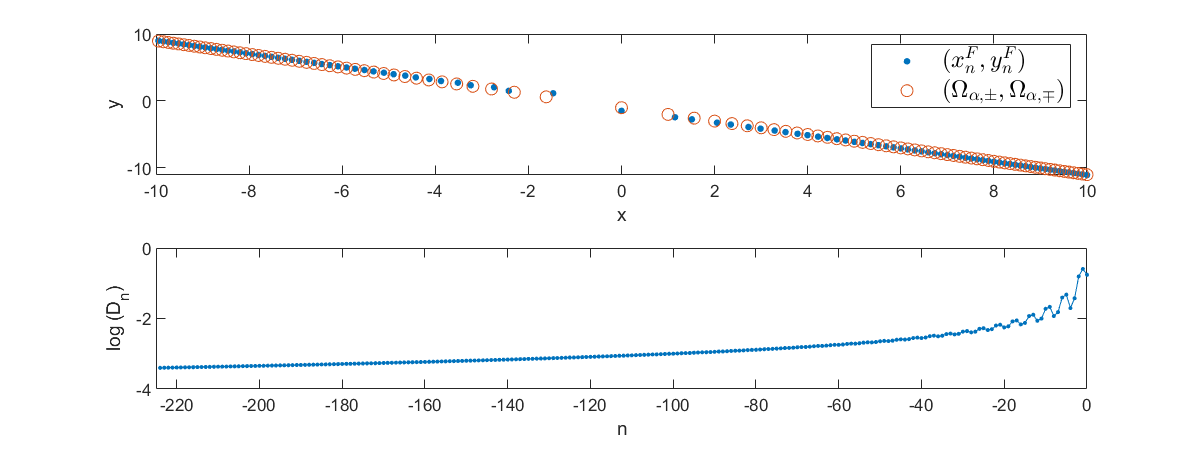}
\caption{Top: Points on the Freud orbit (dots) and in the sequence $\Omega_n$ for $n \le 0$ and $N = 1$. Bottom: Log-distance $\log\left(D_n \right)$ as  a function of $\alpha = n \le 0$.}
\label{fig:p2pts}
\end{figure}

Similar to the previous subsection, we may explicitly determine the period two points of the $\alpha$-dP1 mappings which, in order to be real and genuine, require that $\alpha < 0$. In $(x, y, \alpha)$ coordinates, there is a single period 2 orbit for each value of $\alpha < 0$, given by 
\begin{equation}
(\Omega_{\alpha, \pm},\Omega_{\alpha, \mp}) :=\left(\dfrac{-1 \pm \sqrt{1 - 4 \alpha r }}{2r},\dfrac{-1 \mp \sqrt{1 - 4 \alpha r}}{2r}\right).
\label{eq:alphaper}
\end{equation}
In $(s,f,u)$ coordinates, the period 2 orbit for each $\alpha < 0$ reads
\begin{equation}
(s_\alpha, f_\alpha, u_\alpha) :=\left(0, 0, -2(-1 \pm \sqrt{1 - 4 r \alpha})^{-1}\right),
\label{eq:alphaper2}
\end{equation}
with $u_\alpha$ solving the quadratic equation $\alpha\, r\, u_\alpha^2 - u_\alpha + 1 = 0$. Using $(s,f,u)$ coordinates, we prove in Section \ref{sec:Freud} that orbits that converge to $P_{-\infty}$ along the invariant curve of equation \eqref{eq:CMP2} track the period-2 points given by Equation \eqref{eq:alphaper2} as $n \to -\infty$. In $(x, y, n)$ coordinates, this means that the sequence $x_n$ is asymptotic to the sequence $\Omega_{\alpha,\pm}$ as $n \to -\infty$, where
\begin{equation}
\label{eq:peralpha}
\Omega_{\alpha,\pm}=\left\{\begin{array}{ll}(-1 + \sqrt{1 - 4 \alpha r})/(2r) &\text{ if } n = 2 p\\ & \\ (-1 - \sqrt{1 - 4 \alpha r})/(2r) & \text{ if } n = 2 p + 1 \end{array}\right. \quad \alpha = \frac{n}{N}.
\end{equation}
Similar to Figure \ref{fig:hypfixedpoints}, the top panel of Figure \ref{fig:p2pts} shows points on the Freud orbit for negative values of $n$ (dots), together with the period two points $(\Omega_{\alpha,\pm},\Omega_{\alpha,\mp})$ defined by Equations \eqref{eq:alphaper} and \eqref{eq:peralpha} (open circles). These points alternate between the second and fourth quadrants as $n \to -\infty$. The bottom panel of Figure \ref{fig:p2pts} shows the log of the pointwise distance $D_n$ between the two sequences of points as a function of $n$ for $n \le 0$, where
\[
D_n = \sqrt{(x_n^F - \Omega_{\alpha,\pm})^2 + (y_n^F - \Omega_{\alpha,\mp})^2}, \quad \alpha = \frac{n}{N}.
\]

\subsection{Geometric relationship between autonomous fixed points and Freud's orbit}
\label{sec:actionover}

\begin{figure}[ht]
\centering
\includegraphics[width=0.95\textwidth]{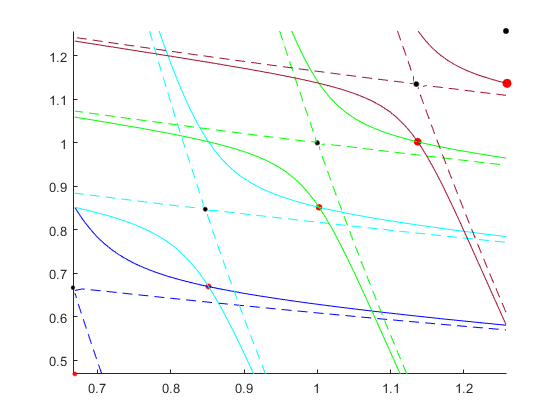}
\caption{$x_{n,N}$ moving in relation to $\omega_{\alpha, +}$ via the intersection and reflection of QRT mappings, for $N=1=r$.}
\label{fig:action}
\end{figure}
As witnessed in the asymptotic expansion \eqref{eq:unP1}, the Freud orbit tracks the fixed points $(w_\alpha,w_\alpha)$ of the autonomous system \eqref{eq:alpha} as $n \rightarrow \infty$. Each of these autonomous mappings is known as a QRT mapping, and their integrable nature provides a geometric explanation for how the Freud orbit and the sequence of points $W_n$ go off to infinity together. We refer the reader to \cite{bib:qrt88} for a thorough discussion of the QRT mappings, but for our purposes it suffices to know that the dynamics of these autonomous systems are restricted to invariant bi-quadratic level sets. Furthermore, the action along these invariant curves has a simple geometric description: to map a point $P$ to its next iterate, one takes $P$'s level set, intersects it with a vertical line through $P$, and then reflects the intersection point over the diagonal. This action can be seen in Figure \ref{fig:action}, where the red point at roughly $(0.85, 0.68)$ maps to the red point at roughly $(1, 0.85)$ by intersecting its level set (the light blue solid curve) with a vertical line (creating an intersection at approximately  $(0.85, 1)$ where the blue and green curves meet) and then reflecting over the diagonal. 

Figure \ref{fig:action} also includes the fixed points {$(w_\alpha,w_\alpha)$ (black dots) and Freud's orbit (red dots) in $(x,y)$ space, both with $n$ ranging from 2 to 6, with $N=r=1$. In this example, when $n=3$, the point $(x_3^F,y_3^F) \approx (0.85, 0.68)$  has a level set whose branches `straddle' the invariant level set of $(w_3,w_3) \approx (0.85, 0.85)$, the dashed separatrix. As discussed above, $(x_3^F,y_3^F)$ maps to $(x_4^F,y_4^F) \approx (1, 0.85)$ via the geometry of intersections and reflections. Now in the non-autonomous dynamics, when  $(x_4^F,y_4^F)$ is to map to  $(x_5^F,y_5^F)$ the mapping \eqref{eq:alpha} has an update of parameter $\alpha$, as $n$ goes from 3 to 4. With this parameter, the invariant bi-quadratic through $(x_4^F,y_4^F)$ is given by the solid green curve, which bears the same relation to $(w_4,w_4)$ as $(x_3^F,y_3^F)$'s level set did to  $(w_3,w_3)$. This process repeats itself in a self-similar fashion, to the invariant red-curves for the next $n$, and so on, ad infinitum.

\section{Numerical simulations}
\label{app:Num}

\begin{figure}[htbp]
\centering
\includegraphics[width=.9\textwidth]{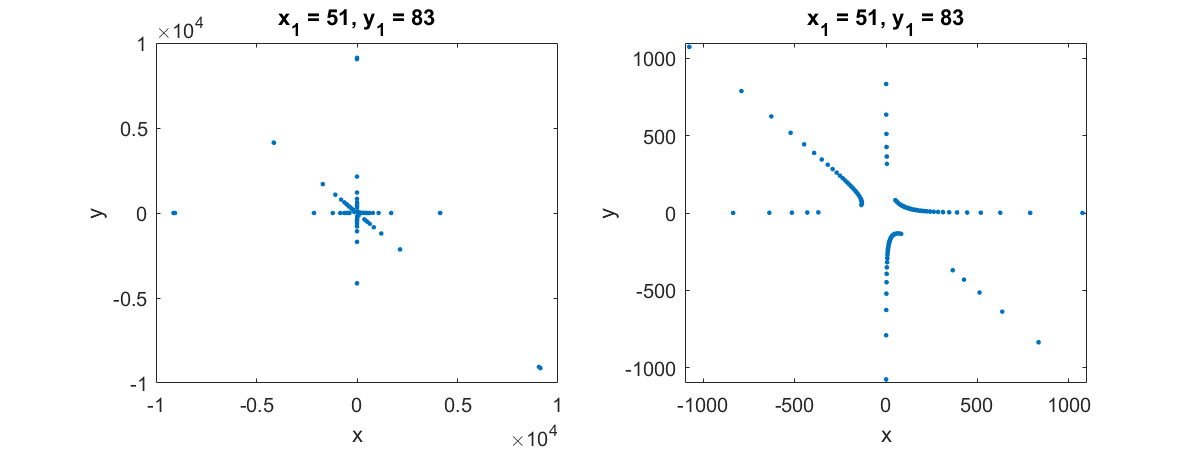}\\
\includegraphics[width=.9\linewidth]{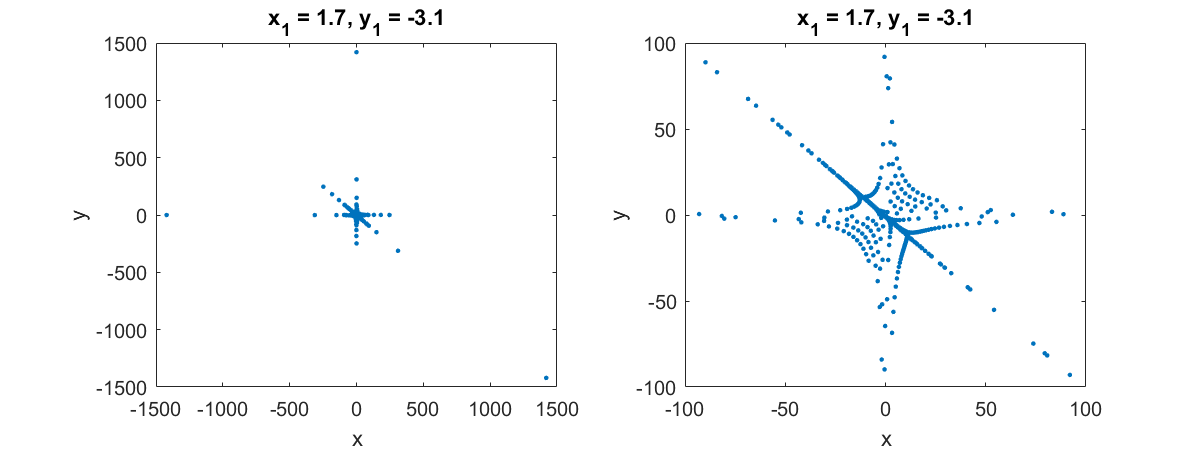}\\
\includegraphics[width=.9\linewidth]{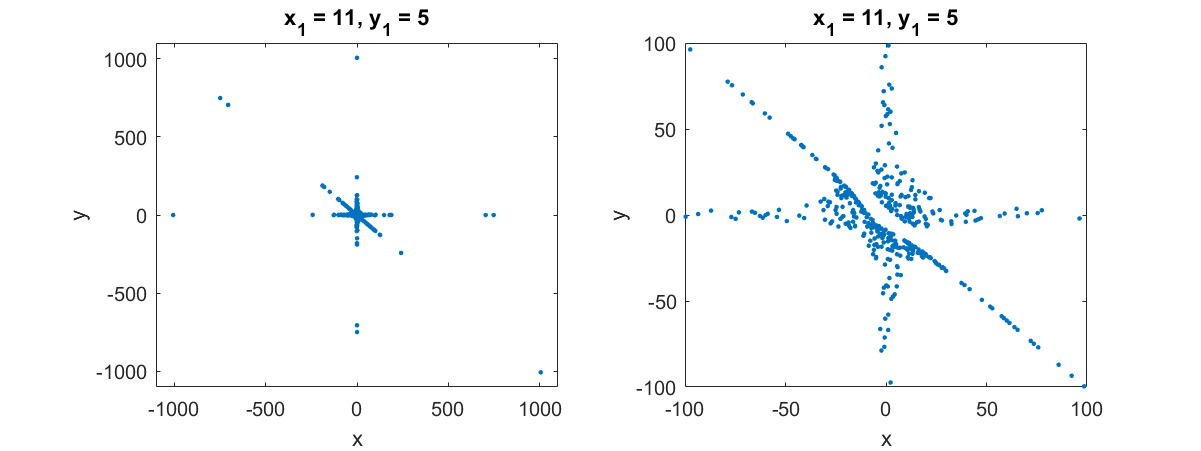}
\caption{Generic orbits of the non-autonomous dP1 equation \eqref{eq:napdintrotwo} with $N = r = 1$. Orbits are plotted in the $(x, y)$ plane, for values of $n \in \{-225, \cdots, 225\}$. The left panels show all of the iterates and the right panel enlarges the region near the origin. Initial condition are $(x_1, y_1) = (51, 83)$ (top), $(x_1, y_1) = (1.7, -3.1)$ (middle), and $(x_1, y_1) = (11, 5)$ (bottom).}
\label{fig:gen_orb}
\end{figure}

Numerical simulations throughout this paper were generated using Python 3.8.5, all using parameters $N=1$ and $r=1$. Note that these parameters can be scaled according to the rescaling presented in \cite{bib:tip20} to provide solutions to an infinite ``ray" of parameters in the $N,r >0$ quadrant. As witnessed in Section \ref{sec:Freud}, the Freud orbit and other orbits with initial conditions on $\mathcal{S}$ do not behave like generic polar orbits, some of which are plotted in Figure \ref{fig:gen_orb}. When such non-polar orbits are computed, numerical error compounds as $n$ gets large, due to accumulation of error along the unstable direction near $P_\infty$. This error grows until these orbits take on negative values for some $x_n$, which is not possible given Theorem \ref{th:LQ} and is thus a numerical artifact. To combat the numerical error, we use the Python library MPMATH, which enables us to perform arithmetic calculations at any desired precision. By computing iterates using thousands of digits of precision, we can reduce the majority of the error to the computation of the initial conditions. Figures are produced in MATLAB, from the Python-generated data.

Initial conditions along  $\mathcal{S}$ were computed using a Lew-Quarles contraction mapping tailored for Equation \eqref{eq:napdintrotwo}, applied to the 0 vector (see Appendix \ref{app:LQ_construction}). With the finite limitations of computation and memory, one cannot initialize using an infinite sequence $\kappa$, but this is not required for a finite number of contractions: for $n$ contractions, one only needs $n+1$ elements in the sequence $\{\kappa_i\}$ to approximate the initial condition $\xi_1$ from a given $\xi_0$. We also note here that, while the contraction mapping could be used to compute later iterates for the same orbit, we instead use the mapping \eqref{eq:napdintrotwo} applied to the initial condition to derive later iterates.

\section{Invariant curves near $P_\infty$ and $P_{-\infty}$}
\label{app:CM}

An approximation of the center manifold $\mathcal C$ near $P_\infty$ may be obtained iteratively by requiring that the curve parametrized by $u$, $s = s_\infty(u), \ f = f_\infty(u)$, remain invariant under the dynamics of \eqref{eq:usf_sys}. Here, $s_\infty$ and $f_\infty$ are polynomial expansions in powers of $u$ and such that the curve $(s_\infty(u), f_\infty(u), u)$ is tangent to the center eigenspace of the linearization of \eqref{eq:usf_sys} near $P_\infty.$ The expressions below are valid to 6th order.
\begin{align}
s_\infty(u)=&2-u-\frac{\gamma}{6} u^2-\frac{\gamma}{36} u^3-\frac{\gamma  (3 \gamma +1)}{216} u^4-\frac{\gamma  (9 \gamma +1)}{1296} u^5 \nonumber \\
&-\frac{\gamma  \left(6 \gamma ^2+18 \gamma +1\right)}{7776} u^6+{\mathcal O}(u^{7}),\nonumber \\
\label{eq:CMP1}\\
f_\infty(u)=&2-u+\frac{\gamma}{6} u^2+\frac{\gamma}{36} u^3-\frac{\gamma  (3 \gamma -1)}{216} u^4-\frac{\gamma  (9 \gamma -1)}{1296} u^5 \nonumber \\
&+\frac{\gamma  \left(6 \gamma ^2-18 \gamma +1\right)}{7776} u^6+{\mathcal O}(u^{7}).
\nonumber
\end{align}

Figure \ref{fig:P1_expansions} compares expansions \eqref{eq:CMP1} to the Freud orbit and provides numerical confirmation that it converges to $P_\infty$ along the center manifold $\mathcal C$.
The above also leads to a parametrization of $\mathcal C$ in $(x,y,n)$ coordinates, which was first obtained in \cite{bib:tip20} and is reproduced below.
\begin{proposition} \label{prop:centerapprox}
The 8\textsuperscript{th} order approximation of the center manifold $\mathcal C$, for $N=1$, parame\-trized in $u$, has the form
\begin{align}
x&=\frac{1}{u}, \nonumber \\
y&=\left[12 r^2 \left(u^2-9 u+3\right) u^5-36 r^3 \left(u^2-9
   u+6\right) u^4+648 r^4 (2-u) u^3-7776 \right.  \nonumber \\
   & \ \ \left.r^5 u^2+46656 r^6-6 r (1-5 u) u^6+u^7\right][46656
   r^6 u]^{-1},  \nonumber \\
   n&=\left[\frac{5 u^7}{3888 r^4}-\frac{u^6}{216r^3}+\frac{u^5}{72
   r^2}-\frac{u^4}{36 r}+3 r+u\right][u^2]^{-1}. \label{eq:centerequation}
\end{align} Numerical estimates of the combinatorial orbit using the center manifold approximation improve as $n$ increases, and as the order of approximation of the manifold increases.
\end{proposition}

\begin{figure}[hbtp]
\centering
\includegraphics[width=.98\textwidth]{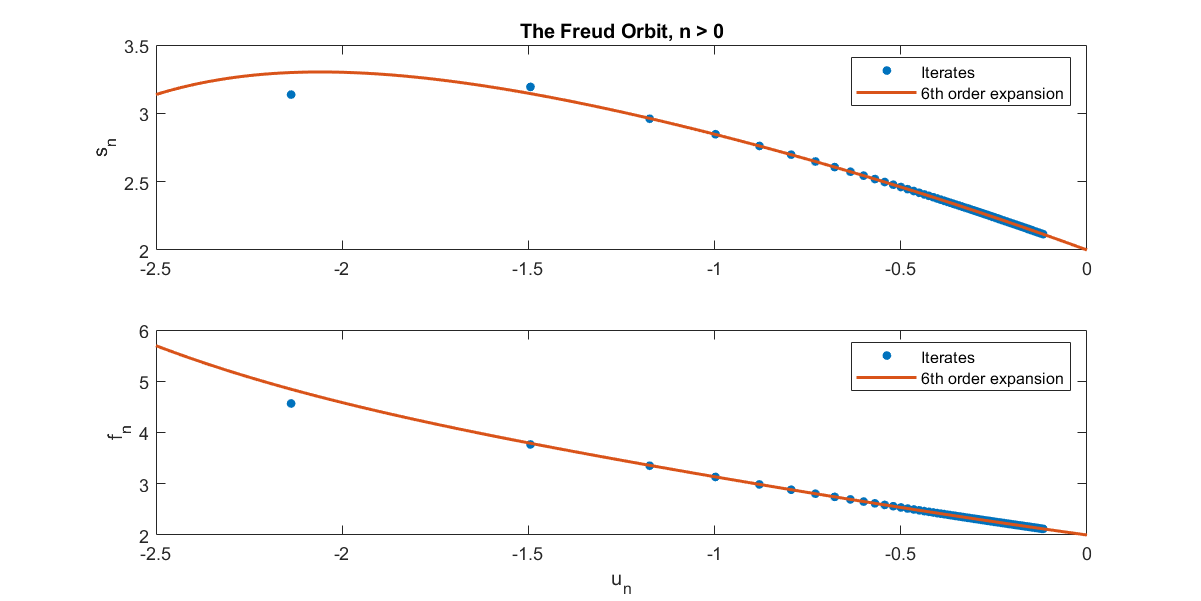}
\caption{Comparison of the numerically obtained iterates of the Freud orbit near $P_\infty$ (dots) with the expansions \eqref{eq:CMP1} (solid curves). The agreement remains quite reasonable even for values of $|u_n|$ of order 1.}
\label{fig:P1_expansions}
\end{figure} 

Similarly, as further explained in the main text, we provide below higher order expressions of the polynomial expansions $s_{-\infty}$ and $f_{-\infty}$ for the invariant curve near $P_{-\infty}$. Corresponding plots are shown in Figure \ref{fig:P2_expansions}, here again providing strong numerical evidence that backward iterates of the Freud orbit converge to $P_{-\infty}$ along the invariant curve $\left(s_{-\infty}(u), f_{-\infty}(u),u\right)$, where
\begin{align}
s_{-\infty}(u)=&-\frac{\gamma}{2} u^2 - \frac{\gamma}{4} u^3 + \frac{\gamma}{8} (\gamma - 1) u^4 + \frac{\gamma}{16} (\gamma - 1) u^5 + \frac{\gamma}{32} (-1 - 2 \gamma + 2 \gamma^2) u^6 \nonumber \\
&+ \frac{\gamma}{64} (-1 - 10 \gamma + 20 \gamma^2) u^7 - \frac{\gamma}{128}(1+25 \gamma -74 \gamma^2+19 \gamma^3)u^8 \nonumber \\
&-\frac{\gamma}{256}(1+49\gamma-168\gamma^2+49\gamma^3) u^9 \nonumber \\
&- \frac{\gamma}{512} (1+84\gamma-252\gamma^2-284\gamma^3+138\gamma^4) u^{10} + {\mathcal O}(u^{11}),\nonumber \\
& \qquad \label{eq:CMP2} \\
f_{-\infty}(u)=&\frac{\gamma}{2} u^2 + \frac{\gamma}{4} u^3 + \frac{\gamma}{8} (\gamma + 1) u^4 + \frac{\gamma}{16} (\gamma + 1) u^5 + \frac{\gamma}{32} (1 - 2 \gamma - 2 \gamma^2) u^6 \nonumber \\
&+ \frac{\gamma}{64} (1 - 10 \gamma - 20 \gamma^2) u^7 - \frac{\gamma}{128}(-1+25 \gamma +74 \gamma^2+19 \gamma^3)u^8 \nonumber \\ 
&-\frac{\gamma}{256}(-1+49\gamma+168\gamma^2+49\gamma^3) u^9\nonumber  \\
&+ \frac{\gamma}{512} (1-84\gamma-252\gamma^2+284\gamma^3+138\gamma^4) u^{10} + {\mathcal O}(u^{11}). \nonumber
\end{align}

\begin{figure}[hbtp]
\centering
\includegraphics[width=.98\textwidth]{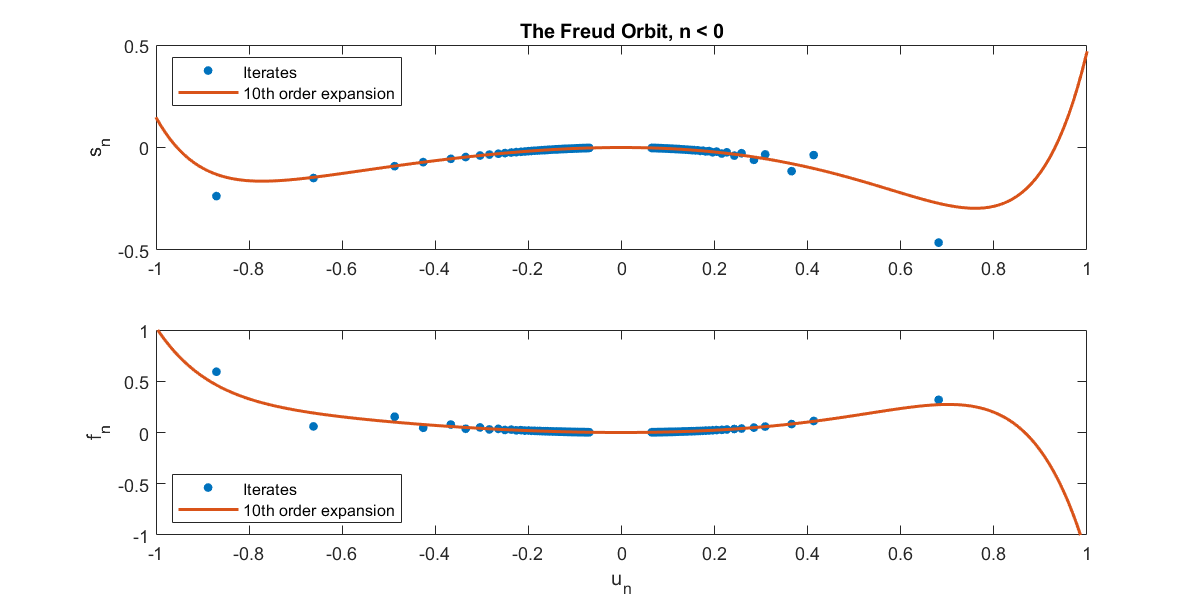}
\caption{Comparison of the numerically obtained backward iterates of the Freud orbit near $P_{-\infty}$ (dots) with the expansions \eqref{eq:CMP2} (solid curves). The agreement remains quite reasonable even for values of $|u_n|$ of order 1.}
\label{fig:P2_expansions}
\end{figure} 
\newpage
\bibliographystyle{siamplain}

\end{document}